\documentclass[12pt]{article}

\usepackage[utf8]{inputenc}
\usepackage[T1]{fontenc}
\usepackage{lmodern}
\usepackage{amsfonts,amsmath,amsthm}

\renewenvironment{proof}{\noindent \textit{Proof.}\newline}{}
\newcommand{\eop}{\hfill $\sqcap\!\!\!\!\sqcup$} 

\usepackage{enumerate}
\usepackage[shortlabels]{enumitem}
\usepackage{color}
\usepackage[all]{xy}
\def\dd{\displaystyle }
\def\n{\mathbb N}
\def\r{\mathbb R}
\def\cc{\mathbb C}
\def\mm{\medskip }
\usepackage{comment}
\allowdisplaybreaks
\usepackage{chngcntr}
\counterwithout{equation}{section} 
\usepackage{graphicx}
\usepackage{dsfont}
\usepackage{titlesec}
\titleformat{\section}
  {\normalfont\Large\bfseries}{\S\thesection}{1em}{}
 \usepackage{textcomp}
 \def\ghil{\textquotesingle}
\usepackage{indentfirst}
\usepackage{cleveref}
\crefname{section}{§}{§§}
\usepackage{authblk}
\usepackage{url}

\begin{document}
\title{Operators on Spaces of Functions and Measures. Vector Invariant (Fractal) Measures}
\date{}

\author{Ion Chi\c tescu, Loredana Ioana\footnote{Corresponding author}, Radu Miculescu \mbox{ and } Lucian Ni\c t\u a}

\newcommand{\Addresses}{{
  \bigskip
  \footnotesize
\textbf{Authors' addresses}

\mm
  Ion Chi\c tescu, \textsc{Faculty of Mathematics and Computer Science, University of Bucharest, Academiei Str.14, 010014, Bucharest, Romania}\par\nopagebreak
  \textit{E-mail address}, Ion Chi\c tescu: \texttt{ionchitescu@yahoo.com}

  \medskip

  Loredana Ioana (Corresponding author), \textsc{Faculty of Mathematics and Computer Science, University of Bucharest, Academiei Str.14, 010014, Bucharest, Romania}\par\nopagebreak
  \textit{E-mail address}, Loredana Ioana: \texttt{loredana.madalina.ioana@gmail.com}

  \medskip

  Radu Miculescu, \textsc{Faculty of Mathematics and Computer Science, University of Bucharest, Academiei Str.14, 010014, Bucharest, Romania}\par\nopagebreak
  \textit{E-mail address}, Radu Miculescu: \texttt{miculesc@yahoo.com}
  
  \medskip
  
  Lucian Ni\c t\u a, \textsc{Technical University of Civil Engineering, Lacul Tei Blvd.,122-124, 020396, Bucharest, Romania}\par\nopagebreak
  \textit{E-mail address}, Lucian Ni\c t\u a: \texttt{luci6691@yahoo.com}

}}







\maketitle

\begin{abstract}
We consider a general schema involving measure spaces, contractions and linear and continuous operators. Within the framework of this schema we use our sesquilinear uniform integral and introduce some integral operators on continuous vector functions spaces, which lead us to operators on spaces of vector measures. Using these last operators, we generalize the Markov operators, obtaining via contractions vector invariant (fractal) measures. Concrete examples are provided.


\end{abstract}

\section{Introduction}

We consider the \underline{fractals} introduced in the spirit of the seminal paper \cite{H2} of J. Hutchinson. The basic theory is exposed in the standard monograph \cite{B} of M.F. Barnsley. See also \cite{F}.

Hutchinson's schema uses \underline{iterated function systems} (IFS), which are finite sets of contractions of a complete metric space $X$. One considers the Hausdorff metric on the complete metric space $\mathcal{K}(X)$ of the non empty compact sets of $X$. A special contraction (the Hutchinson contraction) is constructed, on $\mathcal{K}(X)$, using the aforementioned IFS. The fixed point of the Hutchinson contraction is a \underline{fractal} (in a large sense), e.g. the Cantor set is such a fractal. The Hutchinson schema refers to invariant (fractal) measures, too. Namely, using a probability distribution and the aforementioned IFS, one constructs the so called \underline{Markov operator} acting over the probabilities on the Borel sets of $X$. With a special metric (the Monge-Kantorovich metric), this Markov operator is a contraction on a complete metric space. Its fixed point is an invariant (fractal) measure (probability).

The main goal of this paper is to generalize the theory of invariant (fractal) probabilities for vector measures.

To this end, we developed a preliminary apparatus. The first part of this apparatus appears in detail in our previous papers \cite{CIMN} and \cite{CIMN2}. Namely, in \cite{CIMN} we introduce a sesquilinear uniform integral, which is used in \cite{CIMN2} to define various norms and distances in the space of vector measures of bounded variation.

Our theory is developed within the framework of a general schema involving contractions and operators on Hilbert spaces, all connected in a ''measurable manner''. An intermediate step consists in the construction  and study of some special integral operators on spaces of vector valued continuous functions.
Passing to adjoint operators and using the duality between vector valued continuous functions and (dual) vector valued measures, we arrive to construct and study some special operators on spaces of vector valued measures. These operators viewed within the framework of our general schema lead to  substantial generalizations of the Markov operators from the classic case. Using the distances introduced in the preliminary part on the spaces of vector measures where these Markov-type operators act, we construct  contractions and their fixed points, which are invariant (fractal) vector measures.

The idea of generating invariant (fractal) vector measures (instead of probabilities) appears, in a different form, in some other papers. We quote here some of them, beginning with \cite{MV2}, where a direct generalization of Markov operators for vector measures appears. The reader can also consult \cite{BL}, \cite{FN}, \cite{TM} and \cite{MV}. It is to be mentioned that the discrete case (see \cref{section5}) is more close to the ideas in the aforementioned papers, being a direct generalization of the classic case. 
See also our paper \cite{CIMN3}.

We believe that the study of Markov-type operators on vector measures and of their fixed points has not only purely theoretical reasons. For instance, behaviour of fluids or of electric and magnetic fields can be better described using vector measures.

Let us pass to a brief survey of the contents of the paper. 

In the second paragraph (''Preliminary Facts'') we introduce the notions and results which will be used throughout the paper. A special attention is given to the sesquilinear uniform integral introduced in \cite{CIMN} and to the generalizations of the Monge-Kantorovich and modified Monge-Kantorovich norms and distances introduced in \cite{CIMN2}, making the paper to be self-contained.

In the third paragraph (''Framework of the paper'') we introduce the general schema (framework) which will be used in the sequel. This schema contains a measure space $(\Theta,\Sigma,W)$, a compact metric space $T$, a Hilbert space $X$, Lipschitz functions on $T$, operators on $X$, all measurably indexed with $\theta\in\Theta$. This schema is a general abstractization of many standard models.

In the fourth paragraph (''Special operators on spaces of continuous functions and on spaces of measures'') we use the Bochner integral, the ses\-qui\-li\-near uniform integral and the generalized Monge-Kantorovich (or modified Monge-Kantorovich) norms and distances. First, we introduce some integral operators on vector continuous functions. Using these operators, we pass to adjoints and obtain some operators on spaces of vector measures.

The fifth paragraph (''Particular cases'') studies various cases when the general schema can be applied: the case when all the Lipschitz functions (contractions) are constant, operator semigroups, the discrete case (what concerns $(\Theta,\Sigma,W)$).

The final (sixth) paragraph (''Invariant (fractal) measures'') uses the preceding results concerning operators on spaces of vector measures and the contraction principle to construct invariant (fractal) vector measures, which are fixed points of some Markov-type operators. Concrete examples, together with numerical computations appear.

\section{Preliminary Facts}

Throughout this paper: $\n=\{0,1,2,\dots,n,\dots\},\;\n^*=\{1,2,\dots,n,\dots\},$ $\r_+=[0,\infty)$ and $K$ will be the scalar field (real if $K=\r$, or complex if $K=\cc$). All the sequences will be indexed by $\n$ or $\n^*$ and all the vector spaces (which are assumed to be non null) will be over $K$. We shall write for a sequence $(x_n)_n$ and a non empty set $H:(x_n)_n\subset H,$ to denote the fact that $x_n\in H$ for any $n$.

For any set $T,\mathcal{P}(T)$ is the set of all subsets of $T$. If ${A\subset T,\;\varphi_A:T\rightarrow K}$ is the characteristic (indicator) function of $A$. If $T$ is a non empty set, $X$ is a vector space, $\varphi:T\rightarrow K$ and $f:T\rightarrow X,$ we can consider the function $\varphi f:T\rightarrow X$ defined via $(\varphi f)(t)=\varphi(t)f(t)$ for any $t\in T$ (many times, $f$ will be constant).

If $(E,\|\cdot\|)$ and $(F,|\|\cdot\||))$ are normed spaces, we consider the vector space $\mathcal{L}(E,F)=\{V:E\rightarrow F\mid V\mbox{ is linear and continuous}\}$ normed with the operator norm $\|V\|_o=sup\{|\|V(x)\||\mid x\in E,\|x\|\leq 1\}$ (which is even a Banach space if $F$ is a Banach space). In case $E=F$ we write $\mathcal{L}(E)$ instead of $\mathcal{L}(E,E)$. 
We consider the identity operator $I\in\mathcal{L}(E),$ acting via $I(x)=x$ for any $x\in E$. If $F=K,$ we write $E'$ instead of $\mathcal{L}(E,K)\;(E'$ is the dual of $E$). Considering the normed space $(E,\|\cdot\|)$ (many times we write only $E$), we have the weak* topology $\sigma(E',E)$ of $E'$ (given by the family of seminorms $(\pi_x)_{x\in E}$, where $\pi_x(x')=|x'(x)|,x'\in E'$). For any $V\in\mathcal{L}(E,F)$, the adjoint of $V$ is $V':F'\rightarrow E'$ given via $V'(y')=y'\circ V$ for any $y'\in F'$.

The scalar product of two elements $x,y$ in a Hilbert space $X$ will be denoted by $(x,y)$. In case $X=K^n$, we have, for $x=(x_1,x_2,\dots,x_n)$ and $y=(y_1,y_2,\dots,y_n)$, the standard scalar product $\dd(x,y)=\sum_{i=1}^n x_i\overline{y_i}$, generating the euclidian norm $\dd\|x\|=(\sum_{i=1}^n|x_i|^2)^{\frac 12}$.
For a general Hilbert space $X$ with the scalar product $(\cdot,\cdot)$ and for $V\in\mathcal{L}(X)$, the Hilbert adjoint of $V$ is $V^*\in\mathcal{L}(X)$. (hence $(V(x),y)=(x,V^*(y))$ for any $x,y$ in $X$).

For any non empty set $T$ and any normed space $(X,\|\cdot\|)$, we can consider the Banach space
$$B(T,X)=\{f:T\rightarrow X\mid f\mbox{ is bounded}\}$$
equipped with the norm $f\mapsto\|f\|_\infty=sup\{\|f(t)\|\mid t\in T\}$ (the norm of uniform convergence).

We shall work in the particular situation when $(T,d)$ is a compact metric space ($T$ having at least two elements). Then we have $\mathcal{C}(T,X)\subset B(T,X)$, where $\mathcal{C}(T,X)=\{f:T\rightarrow X\mid f\mbox{ is continuous}\}$ is a Banach space when equipped with the induced norm $\|\cdot\|_\infty$. Many times we write only $\mathcal{C}(X)$ (resp. $B(X)$) instead of $\mathcal{C}(T,X)$ (resp. $B(T,X)$).

Let $(T,d)$ and $(X,\rho)$ be two metric spaces, $T$ having at least two elements and let $f:T\rightarrow X$. The Lipschitz constant of $f$ is defined by the formula
$$\|f\|_L=sup\left\{\frac{\rho\big(f(x),f(y)\big)}{d(x,y)}\;\Big|\; x,y\in T,x\neq y\right\}.$$

In case $\|f\|_L<\infty$, we say that $f$ is lipschitzian. In this case, we have $\rho(f(x),f(y))\leq\|f\|_L\cdot d(x,y)$ for any $x$ and $y$ in $T$. The set of all lipschitzian functions $f:T\rightarrow X$ will be denoted by $Lip(T,X)$. In case $X=T,$ we write $Lip(T)$ instead of $Lip(T,X)$. In the particular case when $X$ is a normed space, it follows that $Lip(T,X)$ is a vector space seminormed with the seminorm $f\mapsto\|f\|_L$. In the particular case when $(T,d)$ is a compact metric space and $X$ is a normed space, it follows that $Lip(T,X)\subset \mathcal{C}(T,X)\subset B(T,X)$ and $Lip(T,X)$ is a normed space with the norm $f\rightarrow\|f\|_{BL}\overset{def}{=}\|f\|_\infty+\|f\|_L$. In the same context, we introduce the sets $L_1(X)=\{f\in Lip(T,X)\mid \|f\|_L\leq 1\}$ and $BL_1(X)=\{f\in Lip(T,X)\mid\|f\|_{BL}\leq 1\}$ (clearly $BL_1(X)\subset L_1(X)$).

A function $f\in Lip(T)$ with $\|f\|_L<1$ is called a contraction (with contraction factor $\|f\|_L$). The fundamental theorem of the fixed point theory is:

\textit{The Contraction Principle (Banach-Cacciopolli-Picard)}

Assume that $(T,d)$ is a complete metric space and $f:T\rightarrow T$ is a contraction. Then $f$ has a unique fixed point $x^*\in X,\mbox{ i.e. } f(x^*)=x^*$.

We use standard facts concerning general measure and integral theory, among them the Bochner integral. Let us mention the fact that, if $\mu$ is an arbitrary positive measure, the space $L^2(\mu)$ with standard norm $\|\cdot\|_2$ is a Hilbert space, the scalar product of two elements $\tilde{f}$ and $\tilde{g}$ in $L^2(\mu)$ being $(\tilde{f},\tilde{g})=\dd\int f\overline{g}d\mu,\mbox{ where }f\in\tilde{f}\mbox{ and }g\in\tilde{g}$ are arbitrary representatives.

Passing to vector measures, we consider an arbitrary non empty set $T$, an arbitrary $\sigma$-algebra of sets $\mathcal{B}\subset \mathcal{P}(T)$ and an arbitrary Banach space $X$. For an arbitrary $\sigma$-additive measure $\mu:\mathcal{B}\rightarrow X$, we define its total variation $|\mu|(T)$.

Let us introduce
$$cabv(\mathcal{B},X)=\{\mu:\mathcal{B}\rightarrow X\mid \mu\mbox{ is }\sigma-\mbox{additive and }|\mu|(T)<\infty\}$$
which becomes a Banach space, when equipped with the variational norm $\mu\rightarrow\|\mu\|=|\mu|(T)$.

For any $0<a<\infty$, write $B_a(X)=\{\mu\in cabv(\mathcal{B},X)\mid \|\mu\|\leq a\}$.

In the present paper, we work in the particular case when $(T,d)$ is a compact metric space (and we shall write $\mathcal{B}\overset{def}{=}$ the Borel sets of $T$). Also, we shall write only $cabv(X)$ instead of $cabv(\mathcal{B},X)$. In general, for any topological space $Z$, the Borel sets of $Z$ will be denoted by $\mathcal{B}_Z$. So, we have the inclusion (for topological spaces $Y,Z$): $\mathcal{B}_Y\otimes\mathcal{B}_Z\subset \mathcal{B}_{Y\times Z}$, where $\mathcal{B}_Y\otimes\mathcal{B}_Z$ is the product $\sigma$-algebra of $\mathcal{B}_Y$ and $\mathcal{B}_Z$ and $Y\times Z$ is the product topological space of $Y$ and $Z$. In case $Y(\mbox{ or }Z)$ is metrizable and separable, one has $\mathcal{B}_Y\otimes\mathcal{B}_Z= \mathcal{B}_{Y\times Z}$.

We continue introducing the basic facts from our previous papers \cite{CIMN} and \cite{CIMN2}. Again $(T,d)$ is a compact metric space and $X$ is a Hilbert space.

A function $f$ of the form $f=\dd\sum_{i=1}^m\varphi_{A_i}x_i$, with $(A_i)_{1\leq i\leq m}$ forming a partition of $T$ and all $x_i\in X$, is called simple. A function $g:T\rightarrow X$ having the property that there exists a sequence $(f_n)_n$ of simple functions such that $f_n\xrightarrow[n]{u}f$ (i.e. $(f_n)_n$ converges uniformly to $f$) is called totally measurable. The vector space of totally measurable functions will be denoted by $TM(X)$. We have the inclusion $\mathcal{C}(X)\subset TM(X)\subset B(X)$.

For any simple function $f=\dd\sum_{i=1}^m\varphi_{A_i}x_i$ and any $\mu\in cabv(X)$, the integral of $f$ with respect to $\mu$ is defined via 
$$\int fd\mu\overset{def}{=}\sum_{i=1}^m (x_i,\mu(A_i)).$$

Then, taking an arbitrary $f\in TM(X)$, we extend the previous definition. Namely, the integral of $f$ with respect to $\mu$ is (coherent definition)
$$\int fd\mu=\underset{m}{lim}\int f_m d\mu,$$
where $(f_m)_m$ is a sequence of simple functions such that $f_m\xrightarrow[m]{u}f$. So our integral is uniform. It is sesquilinear, because the function $(f,\mu)\mapsto\dd\int fd\mu$ is linear in $f$ and antilinear in $\mu$ (when we work with $K=\cc$; for $K=\r$ we have bilinearity). Because of the inequality
$$\Big|\int fd\mu\Big|\leq\|\mu\|\cdot\|f\|_\infty$$
we see that the aforementioned function of $(f,\mu)$ is continuous for $f\in TM(X)$ normed with $\|\cdot\|_\infty$ and $\mu\in cabv(X)$ normed with the variational norm. For any $f\in\mathcal{C}(T,X)$, any $t\in T$ and any $x\in X$, we have $\dd\int fd(\delta_t x)=(f(t),x)$. Here $\delta_t$ is the Dirac measure concentrated at $t$. 

An important interpretation of the integral just introduced is the fact that we have an isometric and antilinear isomorphism (bijection) $H:cabv(X)\rightarrow\mathcal{C}(X)'$ which permits the identification $cabv(X)\equiv\mathcal{C}(X)'$. More details will be given in section B at \S4.

Using this integral, we introduce on $cabv(X)$ and on some of its subspaces new norms (weaker than the variational norm).

For any $\mu\in cabv(X)$, the Monge-Kantorovich norm of $\mu$ is defined via
$$\|\mu\|_{MK}=sup\left\{\Big|\int fd\mu\Big|\mid f\in BL_1(X)\right\}$$
and we get the (generally incomplete) normed space ($cabv(X),\|\cdot\|_{MK}$). For any $\mu\in cabv(X)$ and any $f\in Lip(T,X)$ one has
$$\|\mu\|_{MK}\leq\|\mu\|\mbox{ and }\Big|\int fd\mu\Big|\leq\|\mu\|_{MK}\cdot\|f\|_{BL}.$$

For any $v\in X$, let us define
$$cabv(X,v)=\{\mu\in cabv(X)\mid \mu(T)=v\}.$$

It is clear that $cabv(X,0)$ is a vector subspace of $cabv(X)$ and $\delta_t v\in cabv(X,v)$ for any $t\in T$. It follows that, if $0<a<\infty$ and $v\in X$ is such that $\|v\|\leq a$, then

$$B_a(X,v)\overset{def}{=}B_a(X)\cap cabv(X,v)$$
is not empty, because $\delta_t v\in B_a(X,v)\mbox{ for any }t\in T.$

For any $\mu\in cabv(X,0)$, the modified Monge-Kantorovich norm of $\mu$ is defined via
$$\|\mu\|_{MK}^*\overset{def}{=}sup\left\{\Big|\int fd\mu\Big|\mid f\in L_1(X)\right\}$$
and we get the (generally incomplete) normed space $(cabv(X,0),\|\cdot\|_{MK}^*)$. For any $\mu\in cabv(X,0)$ and any $f\in Lip(T,X)$, one has 
\begin{align*}
\Big|\int fd\mu\Big|&\leq\|\mu\|_{MK}^*\cdot\|f\|_L\\
\|\mu\|_{MK}&\leq\|\mu\|_{MK}^*\leq\|\mu\|_{MK}( diam(T)+1),
\end{align*}
where, as usual, $diam(T)=sup\{d(x,y)\mid x,y\in T\}$.

Using the aforementioned identification $cabv(X)\equiv\mathcal{C}(X)'$, we have the following results, valid for $0<a<\infty,\; n\in\n^*$ and $v\in K^n$ with $\|v\|\leq a$:

The set $B_a(K^n)$, equipped with the metric $d_{MK}$ given via $d_{MK}(\mu,\nu)=\|\mu-\nu\|_{MK}$ and the non empty set $B_a(K^n,v)$, equipped with the metric $d_{MK}$ or with the equivalent metric $d_{MK}^*$ given via $d_{MK}^*(\mu,\nu)=\|\mu-\nu\|_{MK}^*$, are compact metric spaces, their topology being exactly the topology induced by the weak* topology.

In the particular case $K=\r,\;n=1$ and $a=v=1$, the set $B_1^+(\r,1)=B_1(\r,1)\cap\{\mu:\mathcal{B}\rightarrow\r\mid\mu\geq 0\}=$ the probabilities on $\mathcal{B}$, is weak* closed, hence compact for the weak* topology generated by $d_{MK}$ or by $d_{MK}^*$.

For general topology, see \cite{B}. For general measure theory, see \cite{H}. For functional analysis, see \cite{DS}. For vector measures and integration, see \cite{D}.

\section{Framework of the paper}

We shall consider a measure space $(\Theta,\Sigma, W)$ which will be called \underline{the index}\newline\underline{space}, a compact metric space $(T,d)$ (with $card(T)\geq 2$) and a (non null) Hilbert space $X$.

On these spaces, we shall consider the measurable functions $\omega:T\times\Theta\rightarrow X$ and $R:X\times\Theta\rightarrow X$. Namely, $\omega$ is $(\mathcal{B}_T\otimes\Sigma,\mathcal{B}_T)$-measurable and $R$ is $(\mathcal{B}_X\otimes\Sigma,\mathcal{B}_X)$-measurable. We shall use the following indicial notations (for any $\theta\in\Theta):\omega_\theta:T\rightarrow T$ is the function defined via $\omega_\theta(t)=\omega(t,\theta),$ whereas $R_\theta:X\rightarrow X$ is the function defined via $R_\theta(x)=R(x,\theta)$.

We shall assume that, for any $\theta\in\Theta$, one has $R_\theta\in\mathcal{L}(X)$ and $\omega_\theta\in Lip(T),$ with $\|\omega_\theta\|_L=r_\theta$ (in case $r_\theta<1,\;\omega_\theta$ is a contraction)

Before passing further, we shall remark \underline{two particular cases}.

The particular case when \underline{all} $\omega_{\theta},\;\theta\in\Theta,$  \underline{are constant}.

In this case, write $\omega_\theta(t)=t_\theta\in T$ for any $t\in T$ and define $\varphi:\Theta\rightarrow T$ via $\varphi(\theta)=t_\theta$. The fact that $\omega$ is $(\mathcal{B}_T\otimes\Sigma,\mathcal{B}_T)$- measurable  is equivalent to the fact that $\varphi$ is $(\Sigma,\mathcal{B}_T)$-measurable, due to the equality $\omega^{-1}(B)=T\times\varphi^{-1}(B),$ for any $B\in\mathcal{B}_T$.

The other particular case we have in mind is the \underline{discrete case}, when the measure space $(\Theta,\Sigma,W)$ is \underline{discrete}, i.e. either $\Theta=\{1,2,\dots,M\}$ for some $M\in\n^*$, or $\Theta=\n^*$, and, in both cases, $\Sigma=\mathcal{P}(\Theta),W=$ the cardinal measure: $W= card,$ defined via $card(A)=$ the number of elements in A (if $A$ is finite) or $card(A)=\infty$ (if $A$ is infinite).

In this case, to say that $\omega$ is $(\mathcal{B}_T\otimes\Sigma,\mathcal{B}_T)$-measurable means to say that the function $\omega_\theta:T\rightarrow T$ is $(\mathcal{B}_T,\mathcal{B}_T)$-measurable for any $\theta\in\Theta$. Indeed, the measurability of $\omega$ implies the measurability of any $\omega_\theta$. Conversely, if all $\omega_\theta$ are $(\mathcal{B}_T,\mathcal{B}_T)$-measurable, we have for any $B\in\mathcal{B}_T$:
\begin{align*}
\omega^{-1}(B)&=\{(t,\theta)\in T\times\Theta\mid\omega(t,\theta)=\omega_\theta(t)\in B\}=\\
&=\bigcup_{\theta\in\Theta}\Big(\omega_\theta^{-1}(B)\times\{\theta\}\Big)\in\mathcal{B}_T\otimes\mathcal{P}(\Theta) \mbox{ }a.s.o.
\end{align*}

Returning to the general case, we introduce the function $Ind:\Theta\rightarrow\r_+,$ given via $Ind(\theta)=\|R_\theta\|_o$. We shall assume that $Ind$ is $(\Sigma,\mathcal{B}_{\r_+})$-measurable.

Before passing further, let us notice that, in case $X$ is separable, Ind is automatically 
$(\Sigma,\mathcal{B}_{\r_+})$-measurable.
Indeed, if $A\subset X$ is at most countable and dense in $X$, we have for any $a\in A$ the $(\Sigma,\mathcal{B}_X)$-measurable function $\theta\mapsto R_\theta(a)$, hence we have the $(\Sigma,\mathcal{B}_{\r_+})$-measurable function $\theta\mapsto\|R_\theta(a)\|$. Then the function Ind, i.e. the function $\theta\mapsto \|R_\theta\|_o$ is obtained as follows:
\begin{align*}
\|R_\theta\|_o&=sup\{\|R_\theta(x)\|\mid x\in X,\|x\|\leq 1\}=\\
&=sup\{\|R_\theta(a)\|\mid a\in A,\|a\|\leq 1\}.
\end{align*}
The last at most countable supremum is $(\Sigma,\mathcal{B}_{\r_+})-$measurable. We also introduce the function $Lip:\Theta\rightarrow\r_+$, via $Lip(\theta)=r_\theta$.

\textbf{Lemma 3.1}
\textit{
The function $Lip:\Theta\rightarrow\r_+$ is $(\Sigma,\mathcal{B}_{\r_+})$-measurable. 
}

\proof
Assume $T_0\subset T$ is at most countable and dense in $T$.
Write $T_0$ as a ''sequence'': $T_0=\{t_i\mid i\in M\}$ where $\emptyset\neq M\subset\n$ is at most countable and $t_i\neq t_j,\mbox{ if }i\neq j$.

For any $i\neq j$ in $M$, we define $f_{ij}:\Theta\rightarrow\r_+$ via

$$f_{ij}(\theta)=\frac{d\big(\omega_\theta(t_i),\omega_\theta(t_j)\big)}{d(t_i,t_j)}=\frac{d\big(\omega(t_i,\theta),\omega(t_j,\theta)\big)}{d(t_i,t_j)}.$$

All the functions $f_{ij}$ are $(\Sigma,\mathcal{B}_{\r_+})$-measurable. This is seen as follows: the function $\theta\mapsto\big(\omega(t_i,\theta),\omega(t_j,\theta)\big)$ is $(\Sigma,\mathcal{B}_T\otimes \mathcal{B}_T)$-measurable, hence $(\Sigma,\mathcal{B}_{T\times T})$-measurable, because $\mathcal{B}_T\otimes \mathcal{B}_T=\mathcal{B}_{T\times T}$; the function $\big(\omega(t_i,\theta),\omega(t_j,\theta)\big)\mapsto d\big(\omega(t_i,\theta),\omega(t_j,\theta)\big)$ is $(\mathcal{B}_{T\times T},\mathcal{B}_{\r_+})$-measurable, because $d:T\times T\rightarrow\r+$ is continuous.

The function $Lip:\Theta\rightarrow\r_+$ is obtained via
\begin{align*}
Lip(\theta)&=r_\theta=sup\left\{\frac{d\big(\omega_\theta(t),\omega_\theta(s)\big)}{d(t,s)}\;\Big| (t,s)\in T\times T, t\neq s\right\}=\\
&=sup\left\{\frac{d\big(\omega_\theta(t_i),\omega_\theta(t_j)\big)}{d(t_i,t_j)}\;\Big| (i,j)\in (M\times M)\setminus\big\{(i,i)\mid i\in\n\big\}\right\}
\end{align*}
and the last at most countable supremum is 
$(\Sigma,\mathcal{B}_{\r_+})$-measurable.\eop

\mm
The final condition we shall impose is the following:
$$\int_\Theta Ind(\theta)(1+Lip(\theta))dW(\theta)<\infty$$
i.e.
$$\int_\Theta \|R_\theta\|_o\cdot(1+r_\theta)dW(\theta)<\infty.$$

Notice that in case $Lip$ is bounded (i.e. $\underset{\theta\in\Theta}{sup}\;r_\theta<\infty$, which is in particular true if all $\omega_\theta$ are contractions), the last condition means

$$\int_\Theta Ind(\theta)dW(\theta)=\int_\Theta \|R_\theta\|_odW(\theta)<\infty.$$

\section{Special operators on spaces of continuous functions and on spaces of measures}

\indent A. This section is dedicated to special operators on spaces of continuous vector functions.

The construction will be carried on step by step as follows.

a) We show that, for any $f\in\mathcal{C}(T,X)$ and any $t\in T$, the function $U:\Theta\rightarrow X$, given via
$$U(\theta)=(R_\theta\circ f\circ\omega_\theta)(t)$$
is Bochner integrable with respect to $W$.

\textit{Proof of this fact}

Because $\omega$ is $(\mathcal{B}_T\otimes\Sigma,\mathcal{B}_T)$-measurable, it follows that the function $(\tau,\theta)\mapsto f\big(\omega(\tau,\theta)\big)$ is $(\mathcal{B}_T\otimes\Sigma,\mathcal{B}_X)$-measurable, i.e. the function $(\tau,\theta)\mapsto (f\circ\omega_\theta)(t)$ is $(\mathcal{B}_T\otimes\Sigma,\mathcal{B}_X)$-measurable. Consequently, the function $(\tau,\theta)\mapsto\big((f\circ\omega_\theta)(\tau),\theta\big)$ is $(\mathcal{B}_T\otimes\Sigma,\mathcal{B}_X\otimes\Sigma)$-measurable. It follows that the function
\newline
$(\tau,\theta)\mapsto R\big((f\circ\omega_\theta)(\tau),\theta\big)$ is $(\mathcal{B}_T\otimes\Sigma,\mathcal{B}_X)$-measurable. In other words, the function 
$(\tau,\theta)\mapsto R_\theta\big((f\circ\omega_\theta)(\tau)\big)=(R_\theta\circ f\circ\omega_\theta)(\tau)$ is 
$(\mathcal{B}_T\otimes\Sigma,\mathcal{B}_X)$-measurable.

As a consequence, for the already fixed $t\in T$, the function 
\newline
$\theta\mapsto (R_\theta\circ f\circ\omega_\theta)(t)$ is $(\Sigma,\mathcal{B}_X)$-measurable.

We proved that $U$ is $(\Sigma,\mathcal{B}_X)$-measurable.

To finish the proof of the Bochner integrability of $U$, we notice that, for any $\theta\in\Theta$, one has $\|U(\theta)\|=\|(R_\theta\circ f\circ\omega_\theta)(t)\|\leq\|R_\theta\|_o\cdot\|f\|_\infty$, hence 
\begin{align*}
\dd\int\|(R_\theta\circ f\circ\omega_\theta)(t)\|dW(\theta)&\leq\int\|R_\theta\|_o\cdot\|f\|_\infty dW(\theta)\leq\\
&\leq \int\|R_\theta\|_o(1+r_\theta) dW(\theta)<\infty.
\end{align*}

b) The preceding fact enables us to construct, for any $f\in\mathcal{C}(T,X)$, the function $H(f):T\rightarrow X$, given via
$$H(f)(t)=\int(R_\theta\circ f\circ\omega_\theta)(t) dW(\theta).$$

Our next result is that, for any $f\in\mathcal{C}(T,X)$, one has $H(f)\in\mathcal{C}(T,X)$.

\textit{Proof of this fact}

We shall fix an arbitrary $t\in T$ and we shall prove that $H(f)(t_n)\xrightarrow[n]{}H(f)(t)$ for any sequence $(t_n)_n\subset T$ such that $t_n\xrightarrow[n]{}t$.

Due to continuity of $f$ and $\omega_\theta$, we have the pointwise convergence $u_n\xrightarrow[n]{}U$, where $U:\Theta\rightarrow X, U(\theta)=(R_\theta\circ f\circ\omega_\theta)(t)$ and $u_n:\Theta\rightarrow X,u_n(\theta)=(R_\theta\circ f\circ\omega_\theta)(t_n)$. As we just proved, the functions $U$ and $u_n$ are Bochner integrable with respect to $W$ and, for any $\theta$ and $n$:
$$\|u_n(\theta)\|\leq\|R_\theta\|_o\cdot\|f\|_\infty,\;\|U(\theta)\|\leq\|R_\theta\|_o\cdot\|f\|_\infty.$$

The function $\theta\mapsto\|R_\theta\|_o\cdot\|f\|_\infty$ is $W$-integrable and Lebesgue's dominated convergence theorem says that
$$\int u_n(\theta)dW(\theta)\xrightarrow[n]{}\int U(\theta)dW(\theta)\Leftrightarrow H(f)(t_n)\xrightarrow[n]{}H(f)(t).$$

c) The previous result enables us to consider the operator (which is obviously linear) $H_c:\mathcal{C}(T,X)\rightarrow\mathcal{C}(T,X)$, given via
$$H_C(f)=H(f).$$

For any $t\in T$, one has
\begin{align*}
\|H_c(f)(t)\|&=\Big\|\int(R_\theta\circ f\circ\omega_\theta)(t) dW(\theta)\Big\|\leq\\
&\leq\int\|(R_\theta\circ f\circ\omega_\theta)(t)\|dW(\theta)\leq\int\|R_\theta\|_o\cdot\|f\|_\infty dW(\theta),
\end{align*}

hence, for any $f\in\mathcal{C}(T,X)$, one has
\begin{equation}
\tag{4.1}\label{eq:4.1}
\big\|H_C(f)\big\|_\infty\leq\int\|R_\theta\|_o\cdot\|f\|_\infty dW(\theta).
\end{equation}

We have proved

\textbf{Theorem 4.1}
\textit{The operator $H_C:\mathcal{C}(T,X)\rightarrow\mathcal{C}(T,X)$ is linear and continuous with
$$\|H_C\|_o\leq\int\|R_\theta\|_o dW(\theta).$$
}

Now, we restrain the action of $H_C$, letting $H_C$ act only on lipschitzian functions.

Again we proceed step by step.

a) For any $f\in Lip(T,X)$ and any $s,t$ in $T$, one has

$$\|H(f)(s)-H(f)(t)\|\leq\Big(\int\|R_\theta\|_o\cdot r_\theta\; dW(\theta)\Big)\cdot\|f\|_L\cdot d(s,t).$$

\textit{Proof of this fact}

\begin{align*}
\|H(f)(s)-H(f)(t)\|&=\Big\|\int\big[(R_\theta\circ f\circ\omega_\theta)(s)-(R_\theta\circ f\circ\omega_\theta)(t)\big]dW(\theta)\Big\|=\\
&=\Big\|\int R_\theta\Big(f\big(\omega_\theta(s)\big)-f\big(\omega_\theta(t)\big)\Big)dW(\theta)\Big\|\leq\\
&\leq\int\|R_\theta\|_o\cdot\|f\big(\omega_\theta(s)\big)-f\big(\omega_\theta(t)\big)\|dW(\theta)\leq\\
&\leq\int\|R_\theta\|_o\cdot\|f\|_L\cdot  d\big(\omega_\theta(s),\omega_\theta(t)\big)dW(\theta)\leq\\
&\leq\|f\|_L\cdot\int\|R_\theta\|_o\cdot r_\theta\cdot d(s,t)dW(\theta)=\\
&=\|f\|_L\cdot d(s,t)\cdot\int\|R_\theta\|_o\cdot r_\theta\;dW(\theta).
\end{align*}

b) We proved that, for any $f\in Lip(T,X)$, one has $H(f)\in Lip(T,X)$ and

\begin{equation}
\tag{4.2}\label{eq:4.2}
\big\|H(f)\big\|_L\leq\|f\|_L\cdot\int\|R_\theta\|_o\cdot r_\theta\;dW(\theta).
\end{equation}

Consequently, we can introduce the (obviously linear) operator $H_L:Lip(T,X)\rightarrow Lip(T,X)$, given via
$$H_L(f)=H(f).$$

Considering on $Lip(T,X)$ the norm $\|\cdot\|_{BL}$, we have, for any $f\in Lip(T,X)$ (use (\ref{eq:4.1}) and (\ref{eq:4.2})):

\begin{align}
\|H(f)\|_{BL}&=\|H(f)\|_\infty+\|H(f)\|_L\leq\nonumber\\
&\leq\Big(\int\|R_\theta\|_o dW(\theta)\Big)\cdot\|f\|_\infty+\Big(\int\|R_\theta\|_o\cdot r_\theta dW(\theta)\Big)\|f\|_L\leq\nonumber\\
&\leq \Big(\int\|R_\theta\|_o dW(\theta)\Big)\big(\|f\|_\infty+\|f\|_L\big)+\nonumber\\ 
&+
\Big(\int\|R_\theta\|_o r_\theta d W(\theta)\Big)\big(\|f\|_\infty+\|f\|_L\big)=\nonumber\\
&=\Big(\int\|R_\theta\|_o(1+r_\theta)dW(\theta)\Big)\cdot\|f\|_{BL}.\tag{4.3}\label{eq:4.3}
\end{align}

Taking this into account, we proved

\textbf{Theorem 4.2}
\textit{The operator ${H_L:\big(Lip(T,X),\|\cdot\|_{BL}\big)\rightarrow \big(Lip(T,X),\|\cdot\|_{BL}\big)}$ is linear and continuous with
$$\|H_L\|_o\leq\int\|R_\theta\|_o(1+r_\theta)dW(\theta).$$
}

B. This section is dedicated to special operators on spaces of vector measures. We shall use the results of the preceding section.

Before passing further, we introduce some precise notations.

It is a classical result now that the spaces $\mathcal{C}(T,X)'$ and $cabv(T,X')$ are li\-ne\-ar\-ly and isometrically isomorphic (see, e.g. \cite{D}), this result being valid for any Banach space $X$ (not only for Hilbert spaces $X$). Namely, one considers, for any $f\in\mathcal{C}(T,X)$ and any $m'\in\mathcal{C}(T,X')$ the \underline{classical linear integral} $\dd\int^* fdm'\in K$.

(Brief recall: for a simple $X$-valued function $\varphi=\dd\sum_{i=1}^n\varphi_{A_i}x_i$ define $\dd\int^*\varphi dm'=\sum_{i=1}^nm'(A_i)(x_i)$ and extend the integral for a continuous function $f\in\mathcal{C}(T,X)$ such that $f=\underset{n}{lim}f_n$ (unform limit, with $f_n$ simple valued functions) via $\dd\int^*fdm'=\underset{n}{lim}\int^*f_ndm'$). 

Then, one has a linear and isometric isomorphism $\Gamma:cabv(T,X')\rightarrow\mathcal{C}(T,X)'$ given as follows: for any $m'\in cabv(T,X'),\;\Gamma(m'):\mathcal{C}(T,X)\rightarrow K$ acts via $\Gamma(m')(f)=\dd\int^*fdm'$, for any $f\in\mathcal{C}(T,X)$.

In order to express this result in terms of our integral, recall first that for any Hilbert space $X$, one has the antilinear and isometric isomorphism $P:X\rightarrow X'$, acting via $P(y)=y'$, where $y'(x)=(x,y)$, for any $x\in X$ (Riesz-Fr\'echet theorem). This gives the antilinear and isometric isomorphism $\Omega:cabv(T,X)\rightarrow cabv(T,X')$, given via $\Omega(m)=P\circ m$ for any $m\in cabv(T,X)$. Then $\Phi=\Gamma\circ\Omega:cabv(T,X)\rightarrow\mathcal{C}(T,X)'$ is an antilinear and isometric isomorphism which identifies $cabv(T,X)$ and $\mathcal{C}(T,X)'$.

Let us see how $\Phi$ works. For any $m\in cabv(T,X)$, the action of $y'=\Phi(m):\mathcal{C}(T,X)\rightarrow K$ on $f\in\mathcal{C}(T,X)$ is given via

$$y'(f)=\int fdm.$$

The last equality is proved first for simple functions (usual trick) and then for continuous functions, passing to uniform limit (see the definition of classical linear integral). Indeed, if $m\in cabv(T,X)$ and $m'=\Omega(m)$, we have for $f=\dd\sum_{i=1}^n\varphi_{A_i}x_i$:
$$\int^* fdm'=\sum_{i=1}^n m'(A_i)(x_i)=\sum_{i=1}^n P\big(m(A_i)\big)(x_i)=\sum_{i=1}^n\big(x_i,m(A_i)\big)=\int fdm.$$

Synthetically, we have the formula
\begin{equation}
\tag{4.4}\label{eq:4.4}
\Phi(m)(f)=\int fdm
\end{equation}
valid for any $m\in cabv(T,X)$ and any $f\in\mathcal{C}(T,X)$.

Taking into account the fact that $\Phi$ is a bijection, we have $m_1=m_2\Leftrightarrow\Phi(m_1)=\Phi(m_2)$, hence (\ref{eq:4.4}) can be rewritten as follows: if $m_1,m_2$ are in $cabv(T,X)$, then 
\begin{equation}
\tag{4.4\ghil}\label{eq:4.4'}
m_1=m_2\Rightarrow\int fdm_1=\int fdm_2\mbox{ for any }f\in\mathcal{C}(T,X).
\end{equation}

Now, let us return to our subject. The operator $H_C:\mathcal{C}(T,X)\rightarrow\mathcal{C}(T,X)$ generates its adjoint $H_C':\mathcal{C}(T,X)'\rightarrow\mathcal{C}(T,X)'$ acting via $H_C'(y)=y'\circ H_C$, for any $y\in\mathcal{C}(T,X)'$.

We introduce the commutative diagram
\begin{displaymath}
\xymatrix{cabv(T,X)\ar[d]_\Phi\ar[r]^{\mathcal{H}}& cabv(T,X)\ar[d]^\Phi\\
\mathcal{C}(T,X)'\ar[r]_{H_C'}&\mathcal{C}(T,X)'
\ar@<1ex>[u]^{\Phi^{-1}}}
\end{displaymath}
where $\mathcal{H}:cabv(T,X)\rightarrow cabv(T,X)$ is defined via
$$\mathcal{H}=\Phi^{-1}\circ H_C'\circ\Phi.$$

Because $\Phi$ and $\Phi^{-1}$ are antilinear, it follows that $\mathcal{H}$ is a linear and continuous operator.

The commutativity of the diagram means
\begin{equation}
\tag{4.5}\label{eq:4.5}
\Phi\circ\mathcal{H}=H_C'\circ\Phi.
\end{equation}

\textbf{Theorem 4.3} (Change of Variable Theorem)
\textit{For any $f\in\mathcal{C}(T,X)$ and any $\nu\in cabv(T,X)$ one has
$$\int fd\mathcal{H}(\nu)=\int H_C(f)d\nu.$$
}

\proof
Let $f\in\mathcal{C}(T,X)$ and $\nu\in cabv(T,X)$. According to (\ref{eq:4.4}) we have $\dd\int fdm=\Phi(m)(f)$ for any $m\in cabv(T,X)$, in particular for $m=\mathcal{H}(\nu)$. Hence (use (\ref{eq:4.5}))
\begin{align*}
\int fd\mathcal{H}(\nu)&=\Phi\big(\mathcal{H}(\nu)\big)(f)=(\Phi\circ\mathcal{H})(\nu)(f)=\\
&=(H_C'\circ\Phi)(\nu)(f)=H_C'\big(\Phi(\nu)\big)(f)=\big(\Phi(\nu)\circ H_C\big)(f)=\\
&=\Phi(\nu)\big(H_C(f)\big)=\int H_C(f) d\nu.\\
& \mbox{ (the final equality with (\ref{eq:4.4})})\hspace{2cm}\Box
\end{align*}

We continue giving some evaluations of the norms of the operator $\mathcal{H}$, viewed as acting in $cabv(T,X)$ or in some subspaces of $cabv(T,X)$, with different norms.
Notations of the type $\|\mathcal{H}\|_{o,norm}$ will be used.

We start naturally with $cabv(T,X)$, equipped with the usual variational norm.

\textbf{Theorem 4.4}
\textit{The operator $\mathcal{H}:\big(cabv(T,X),\|\cdot\|\big)\rightarrow\big(cabv(T,X),\|\cdot\|\big)$ is linear and continuous. We have
$$\|\mathcal{H}\|_{o,var}\leq\int\|R_\theta\|_o dW(\theta).$$
}

\proof
Fix arbitrarily $\nu\in cabv(T,X)$. We have $\|\mathcal{H}(\nu)\|=\|\Phi\big(\mathcal{H}(\nu)\big)\|$.

Then (use (\ref{eq:4.5}), (\ref{eq:4.4}) and (\ref{eq:4.1})):
\begin{align*}
\|\Phi\big(\mathcal{H}(\nu)\big)\|&=\|(H_C'\circ\Phi)(\nu)\|=\|H_C'\big(\Phi(\nu)\big)\|=\\
&=\underset{\|f\|_\infty\leq 1}{sup}\big|H_C'\big(\Phi(\nu)\big)(f)\big|=\underset{\|f\|_\infty\leq 1}{sup}\big|\Phi(\nu)\big(H_C(f)\big)\big|=\\
&=\underset{\|f\|_\infty\leq 1}{sup}\Big|\int H_C(f) d\nu\Big|\leq\underset{\|f\|_\infty\leq 1}{sup}\big\|H_C(f)\big\|_\infty\cdot\|\nu\|\leq\\
&\leq \underset{\|f\|_\infty\leq 1}{sup}\Big(\int\|R_\theta\|_o dW(\theta)\Big)\cdot\|f\|_\infty\cdot\|\nu\|.
\end{align*}

It follows that
\begin{align*}
\|\mathcal{H}\|_{o,var}&=\underset{\|\nu\|\leq 1}{sup}\big\|\mathcal{H}(\nu)\big\|\leq\underset{\|f\|_\infty\leq 1}{sup}\Big(\int \|R_\theta\|_o\cdot dW(\theta)\Big)\cdot\|f\|_\infty\leq\\
&\leq\int\|R_\theta\|_o dW(\theta).
\end{align*}
\eop

\mm
Now, working with the Monge-Kantorovich norm, we obtain

\textbf{Theorem 4.5}
\textit{The operator ${\mathcal{H}:\big(cabv(T,X),\|\cdot\|_{MK}\big)\rightarrow\big(cabv(T,X),\|\cdot\|_{MK}\big)}$ is linear and continuous. We have
$$\|\mathcal{H}\|_{o,MK}\leq\int\|R_\theta\|(1+r_\theta)dW(\theta).$$
}

\proof
Let $\nu\in cabv(T,X)$. We shall use Theorem 4.3 and (\ref{eq:4.3}), obtaining suc\-ce\-ssi\-ve\-ly 
\begin{align*}
\|\mathcal{H}(\nu)\|_{MK}&=\underset{\|f\|_{BL}\leq 1}{sup}\Big|\int fd\mathcal{H}(\nu)\Big|=\underset{\|f\|_{BL}\leq 1}{sup}\Big|\int H_C(f)d\nu\Big|=\\
&=\underset{\|f\|_{BL}\leq 1}{sup}\Big|\int H(f)d\nu\Big|.
\end{align*}

Due to the inequality (valid for any $f\in Lip(T,X)$)
$$\Big|\int H(f)d\nu\Big|=\Big|\int H_L(f)d\nu\Big|\leq\|H_L(f)\|_{BL}\cdot\|\nu\|_{MK},$$
we get
\begin{align*}
\|\mathcal{H}(\nu)\|_{MK}&\leq\|\nu\|_{MK}\cdot\underset{\|f\|_{BL}\leq 1}{sup}\|H(f)\|_{BL}\leq\\
&\leq\|\nu\|_{MK}\cdot\int\|R_\theta\|\cdot(1+r_\theta)\;dW(\theta)\; a.s.o.
\end{align*}
\eop

\mm
In order to use the modified Monge-Kantorovich norm, we need the fo\-llo\-wing intermediary step

\textbf{Lemma 4.6}
\textit{For any $\nu\in cabv(T,X,0)$, one has $\mathcal{H}(\nu)\in cabv(T,X,0).$}

\proof
Take an arbitrary $\nu\in cabv(T,X,0)$. We must prove that $\mathcal{H}(\nu)(T)=0$, i.e. one has $\big(x,\mathcal{H}(\nu)(T)\big)=0$ for any $x\in X$.

To this end, take arbitrarily $x\in X$ and define the constant function $f\in\mathcal{C}(T,X)$, acting via $f(t)=x$ for any $t\in T$. Hence $f=\varphi_T x$ (it is a simple function).

For an arbitrary $t\in T$, one has

$$H_C(f)(t)=\int(R_\theta\circ f\circ\omega_\theta)(t)dW(\theta)=\int R_\theta(x)dW(\theta)\overset{def}{=}y.$$

Hence $H_C(f)$ is a constant function, namely $H_C(f)=\varphi_T y$. With Theorem 4.3, we get

$$\int fd\mathcal{H}(\nu)=\int H_C(f)d\nu\Leftrightarrow\int\varphi_T xd\mathcal{H}(\nu)=\int\varphi_T yd\nu$$
which means
$$\big(x,\mathcal{H}(\nu)(T)\big)=\big(y,\nu(T)\big)$$
and this implies $\Big(x,\mathcal{H}\big(\nu(T)\big)\Big)=0$, because $\nu(T)=0$.\eop

\mm
This invariance result shows that one can consider the ''compressed'' operator $\mathcal{H}_o:cabv(T,X,0)\rightarrow cabv(T,X,0)$ defined via $\mathcal{H}_o(\nu)=\mathcal{H}(\nu)$, for any $\nu\in cabv(T,X,0)$.

We use the modified Monge-Kantorovich norm for this operator $\mathcal{H}_o$, obtaining

\textbf{Theorem 4.7}
\textit{The operator $${\mathcal{H}:\big(cabv(T,X,0),\|\cdot\|_{MK}^*\big)\rightarrow \big(cabv(T,X,0),\|\cdot\|_{MK}^*\big)}$$ is linear and continuous. We have:
$$\|\mathcal{H}\|_{o,MK^*}\leq\int\|R_\theta\|_o\cdot r_\theta\;dW(\theta).$$
}

\proof
We shall use Theorem 4.3 and (\ref{eq:4.2}), obtaining successively for an arbitrary $\nu\in cabv(T,X,0))$ (hence $\mathcal{H}(\nu)\in cabv(T,X,0)$ with Lemma 4.6):
$$\|\mathcal{H}(\nu)\|_{MK}^*=\underset{\|f\|_L\leq 1}{sup}\Big|\int fd\mathcal{H}(\nu)\Big|=\underset{\|f\|_L\leq 1}{sup}\Big|\int H(f)d\nu\Big|.$$

In view of the inequality (valid for any $f\in Lip(T,X)$)
$$\big|\int H(f)d\nu\Big|\leq\|H(f)\|_L\cdot\|\nu\|_{MK}^*,$$
we get
\begin{align*}
\|H(\nu)\|_{MK}^*&\leq\underset{\|f\|_L\leq 1}{sup}\|f\|_L\cdot\Big(\int\|R_\theta\|\cdot r_\theta dW(\theta)\Big)\cdot\|\nu\|_{MK}^*\leq\\
&\leq\|\nu\|_{MK}^*\cdot\int\|R_\theta\|\cdot r_\theta dW(\theta).\; a.s.o.
\end{align*}
\eop

\section{Particular cases}
\label{section5}
In this paragraph we shall study some particular cases of the results obtained in the 
preceding paragraph.

A. We consider the case when \underline{all functions $\omega_\theta:T\rightarrow T,\theta\in\Theta$ are constant}. We saw that, in this case, one has a $(\Sigma,\mathcal{B}_T)$-measurable function $\varphi:\Theta\rightarrow T$ such that $\omega_\theta(t)=t_\theta=\varphi(\theta)$ for any $\theta\in\Theta$ and any $t\in T$. Because $r_\theta=0$ for any $\theta\in\Theta$, we have also $\int_0^\infty\|R_\theta\|_od\theta<\infty$.

In order to follow the action of $\mathcal{H}$ in this case, we notice first that, for any $f\in\mathcal{C}(T,X)$, any $\theta\in\Theta$ and any $t\in T$, one has

$$H_C(f)(t)=\int(R_\theta\circ f\circ\omega_\theta)(t)dW(\theta)=\int R_\theta\Big(f\big(\varphi(\theta)\big)\Big)dW(\theta)\in X$$
and the function $H(f)$ is constant.

In order to continue, we need the following

\underline{Fact}
For any $V:\Theta\rightarrow X$ which is Bochner integrable with respect to $W$ and for any $x\in X$, one has
$$\Big(\int V(\theta)dW(\theta),x\Big)=\int\big(V(\theta),x\big)dW(\theta).$$
(the left integral is Bochner and the right integral is abstract Lebesgue). This fact is proved in the same way as the equality
$$S\Big(\int V(\theta)dW(\theta)\Big)=\int(S\circ V)(\theta)dW(\theta),$$
valid for any Banach space $Y$ and any $S\in\mathcal{L}(X,Y).$

Returning to the main topics, let $f\in\mathcal{C}(T,X)$ and $\nu\in cabv(T,X)$. With Theorem 4.3:
$$\int fd\mathcal{H}(\nu)=\int H(f)d\nu.$$
But $H(f)=\varphi_T\cdot\dd\int R_\theta\Big(f\big(\varphi(\theta)\big)\Big)dW(\theta)$, hence
\begin{align*}
\int H(f)d\nu&=\Big(\int R_\theta\Big(f\big(\varphi(\theta)\big)\Big)dW(\theta),\nu(T)\Big)=\\
&=\int\Big(R_\theta\Big(f\big(\varphi(\theta)\big)\Big),\nu(T)\Big)dW(\theta)=\int\Big(f\big(\varphi(\theta)\big),R_\theta^*\big(\nu(T)\big)\Big)dW(\theta),
\end{align*}
leading to the final result 
\begin{equation}
\tag{5.1}\label{eq:5.1}
\int fd\mathcal{H}(\nu)=\int\Big(f\big(\varphi(\theta)\big),R_\theta^*\big(\nu(T)\big)\Big)dW(\theta),
\end{equation}
valid for any $f\in\mathcal{C}(T,X)$ and any $\nu\in cabv(T,X)$.

Arguing about relation (\ref{eq:5.1}), one sees that the value of $\dd\int fd\mathcal{H}(\nu)$ depends only upon the value $\nu(T)$.

In view of (\ref{eq:4.4'}), if we consider $\nu_1$ and $\nu_2$ in $cabv(T,X)$, we have
$$\mathcal{H}(\nu_1)=\mathcal{H}(\nu_2)\Leftrightarrow\int fd\mathcal{H}(\nu_1)=\int fd\mathcal{H}(\nu_2)$$
for any $f\in\mathcal{C}(T,X)$. It follows that
$$\nu_1(T)=\nu_2(T)\Rightarrow\mathcal{H}(\nu_1)=\mathcal{H}(\nu_2).$$

B. In this section, we consider \underline{operator semigroups.}

Namely, let $X$ be a Banach space and recall that a \underline{uniformly continuous}
\underline{operator semigroup on $X$} is a function $P:[0,\infty)\rightarrow\mathcal{L}(X)$ having the following properties:
\begin{enumerate}[a),noitemsep,nosep]
\item $P$ is continuous;
\item $P(0)=I$;
\item $P(s+t)=P(s)\circ P(t)$, for any $s,t$ in $[0,\infty)$.
\end{enumerate}
General theory asserts the existence (in $\mathcal{L}(X)$) of the limit
$$\underset{t\to 0}{lim}\frac1t\big(P(t)-P(0)\big)\overset{def}{=}A\in\mathcal{L}(X).$$
We call $A$ \underline{the generator of the semigroup}.

The analogue of the well-known additivity theorem of Cauchy says that, for any $t\in[0,\infty)$, one has
$$P(t)=exp(tA)\overset{def}{=}I+\sum_{n=1}^\infty\frac{1}{n!}(tA)^n$$
(convergence in $\mathcal{L}(X)$).

For instance, if one takes $P(t)=e^{-t}I$, for any $t\in[0,\infty)$, we get $A=-I$, hence
$$P(t)=exp(-tI)\mbox{ and }\|R(t)\|_o=e^{-t}\mbox{for any }t\in[0,\infty).$$

In order to work within our general framework, we consider again a Hilbert space $X$ and a uniformly continuous operator semigroup $P:[0,\infty)\rightarrow\mathcal{L}(X)$ on $X$.

Then take $(\Theta,\Sigma, W)$ as follows: $\Theta=[0,\infty),\Sigma=\mathcal{B}_{[0,\infty)}$ and $W=$ the Lebesgue measure on $[0,\infty)$.

We define $R:X\times[0,\infty)\rightarrow X$ via $R(x,\theta)=P(\theta)(x)$, for any $x\in X$ and any $\theta\in[0,\infty)$. Then we shall write $R_\theta\overset{def}{=}P(\theta)$ for any $\theta\in[0,\infty)$, identifying $(R_\theta)_{\theta\in[0,\infty)}\equiv P$.

It is seen that $R$ is continuous, because, in case $x_n\xrightarrow[]{n}x$ in $X$ and $t_n\xrightarrow[]{n}t$ in $[0,\infty)$, one has $R(x_n,\theta_n)\xrightarrow[]{n}R(x,\theta):$
\begin{align*}
&\|R(x_n,\theta_n)-R(x,\theta)\|=
\|R_{\theta_n}(x_n)-R_\theta(x)\|\leq\|R_{\theta_n}(x_n)-R_{\theta_n}(x)\|+\\
&+\|R_{\theta_n}(x)-R_\theta(x)\|\leq
\|R_{\theta_n}\|_o\cdot\|x_n-x\|+\|R_{\theta_n}-R_\theta\|_o\cdot\|x\|\leq\\
&\leq\big(\|R_\theta\|_o+\delta\big)\cdot\|x_n-x\|+\|R_{\theta_n}-R_\theta\|_o\cdot\|x\|,
\end{align*}
where $\delta$ can be taken arbitrarily small (for $n\geq n(\delta)$ great enough, because $R_{\theta_n}\xrightarrow[]{n}R_\theta$, hence $\|R_{\theta_n}\|_o\xrightarrow[]{n}\|R_\theta\|_o$).

The continuity of $R$ implies the $(\mathcal{B}_X\otimes\mathcal{B}_{[0,\infty)},\mathcal{B}_X)$-measurability of $R$, because $\mathcal{B}_X\otimes\mathcal{B}_{[0,\infty)}=\mathcal{B}_{X\times[0,\infty)}$.

We complete the schema in the framework taking $T=[0,1]$ and defining $\omega:[0,1]\times[0,\infty)\rightarrow[0,1]$ as follows: take an arbitrary lipschitzian function $u:[0,1]\rightarrow[0,1]$, a continuous function $a:[0,\infty)\rightarrow [0,1]$ with $a(0)=1$ and $a(\theta)>0$ for any $\theta\in[0,\infty)$ (e.g. take $a(\theta)=\dd\frac{1}{1+\theta}$) and define $\omega:[0,1]\times[0,\infty)\rightarrow[0,1]$ via $\omega(t,\theta)=a(\theta)\cdot u(t)$.

Then, for any $\theta\in[0,\infty),\omega_\theta:[0,1]\rightarrow[0,1]$ acts via $\omega_\theta(t)=a(\theta)\cdot u(t)$, in particular $\omega_0=u$.

The $(\mathcal{B}_{[0,1]}\times\mathcal{B}_{[0,\infty)},\mathcal{B}_{[0,1]})$-measurability of $\omega$ is due to its continuity and to the equality
$\mathcal{B}_{[0,1]}\otimes\mathcal{B}_{[0,\infty)}=\mathcal{B}_{[0,1]\times[0,\infty)}$.

It is seen that, for any $\theta\in[0,\infty)$, one has $r_\theta=a(\theta)\cdot\|u\|_L$.
This implies that the condition $\dd\int_0^\infty\|R_\theta\|_o(1+r_\theta)dW(\theta)<\infty$ is equivalent to the condition
$$\int_0^\infty\|R_\theta\|_od\theta=\int_0^\infty\|P(\theta)\|_od\theta<\infty.$$

In the particular case when $R_\theta=P(\theta)=e^{-\theta}I$ for any $\theta\in[0,\infty)$, we have $\|R_\theta\|_o=e^{-\theta}$ and the condition is fulfilled.

Working in this particular case, we see that, for any $f\in\mathcal{C}(T,X)=\mathcal{C}([0,1],X)$ and any $t\in[0,1]$, one has
$$H(f)(t)=\int_0^\infty(R_\theta\circ f\circ\omega_\theta)(t)dW(\theta)=\int_0^\infty e^{-\theta}f\big(a(\theta)u(t)\big)dW(\theta)$$
(For instance: if $f(t)=tx$, for any $t\in[0,1]$, where $x\in X$ is fixed, we have
$$H(f)(t)=\big(u(t)\int_0^\infty e^{-\theta}a(\theta)d\theta\big)x\mbox{ ).}$$

Also in the particular case $R_\theta=e^{-\theta}I$, we can consider the situation when $u\equiv 1$ (hence all functions $\omega_\theta$ are constant and $\omega_\theta(t)=a(\theta)=\varphi(\theta)$ for any $t\in[0,1]$). Using formula (\ref{eq:5.1}), we obtain for any $\mathcal{C}([0,1],X)$ and any $\nu\in cabv([0,1],X)$:
\begin{equation}
\tag{5.1\ghil}\label{eq:5.1'}
\int fd\mathcal{H}(\nu)=\int_0^\infty\Big(f\big(a(\theta)\big),e^{-\theta}\nu(T)\Big)d\theta.
\end{equation}

Considering the more particular case when $a$ is constant too, i.e. $a(\theta)=t_0$ (for some $t_0\in T$), it follows that, for any $f\in\mathcal{C}(T,X)$:

\begin{align*}
\int fd\mathcal{H}(\nu)&=\int_0^\infty\big(f(t_0),e^{-\theta}\nu(T)\big)d\theta=\big(f(t_0),\nu(T)\big)\cdot\int_0^\infty e^{-\theta}d\theta=\\
&=\big(f(t_0),\nu(T)\big)=\int fd\big(\delta_{t_0}\nu(T)\big).
\end{align*}

Using (\ref{eq:4.4'}), we get from (\ref{eq:5.1'}) that, for any $\nu\in cabv(T,X)$ one has
\begin{equation}
\tag{5.1\ghil\ghil}\label{eq:5.1''}
\mathcal{H}(\nu)=\delta_{t_0}\nu(T).
\end{equation}

From (\ref{eq:5.1''}) we deduce that, for any $\nu\in cabv(T,X)$, the measure $\delta_{t_0}\nu(T)$ is a fixed point of $\mathcal{H}$, i.e.
$$\mathcal{H}\big(\delta_{t_0}\nu(T)\big)=\delta_{t_0}\nu(T).$$

C. This section is dedicated to the \underline{discrete} case. We shall be able to compute effectively $\mathcal{H}(\nu)$ for a given $\nu$.

C1. \underline{The finite case}

In the general framework schema, we take $\Theta=\{1,2,\cdots,M\}$ for some $M\in\n^*,\Sigma=\mathcal{P}(\Theta)$ and $W=card$.

Recall first the definition of the \underline{transported measure} (adapted for the present situation). Let $V:T\rightarrow T$ be a continuous function and let $\mu\in cabv(T,X)$. Then the \underline{transported measure} $V(\mu):\mathcal{B}_T\rightarrow X$ is given via 
$V(\mu)(A)\overset{def}{=}\mu\big(V^{-1}(A)\big)$, for any $A\in\mathcal{B}_T$. It is seen that $V(\mu)\in cabv(T,X)$, more precise $\|V(\mu)\|\leq\|\mu\|$. 

The last inequality is proved as follows. If $(A_1,A_2,\dots,A_n)$ is a partition of $T (A_i\in \mathcal{B}_T$, disjoint, $\dd\bigcup_{i=1}^n A_i=T$), we get the partition $\big(V^{-1}(A_1),V^{-1}(A_2),\newline\dots,V^{-1}(A_n)\big)$ of $T$ with $V^{-1}(A_i)\in\mathcal{B}_T$. Clearly
$$\sum_{i=1}^n\|V(\mu)(A_i)\|=\sum_{i=1}^n\|\mu\big(V^{-1}(A_i)\big)\|\leq|\mu|(T)=\|\mu\|,$$
hence $\|V(\mu)\|\leq\|\mu\|$.

We shall use the following three facts, valid for any $f\in\mathcal{C}(T,X)$ and any $\mu\in cabv(T,X)$.

\begin{enumerate}[a),noitemsep,nosep]
\item For any continuous $V:T\rightarrow T$, one has
$$\int f d\big(V(\mu)\big)=\int(f\circ V)d\mu.$$
This is easily seen for simple $f=\dd\sum_{i=1}^n\varphi_{A_i}x_i$ (because $f\circ V=\dd\sum_{i=1}^n\varphi_{V^{-1}(A_i)}x_i$) and one passes to uniform limit.
\item For any $R\in\mathcal{L}(X)$ one has $\dd\int(R\circ f)d\mu=\int fd(R^*\circ\mu)$ (same procedure).
\item For any $f\in\mathcal{C}(T,X)$, one has $H(f)(t)=\dd\sum_{i=1}^M R_i\circ f\circ\omega_i$. Indeed, for any $t\in T$,
$$H(f)(t)=\int(R_\theta\circ f\circ\omega_\theta)(t)dW(\theta)=\sum_{i=1}^M(R_i\circ f\circ\omega_i)(t).$$
\end{enumerate}

\textbf{Theorem 5.1}
\textit{In the context from above, one has, for any $\nu\in cabv(T,X)$:
$$\mathcal{H}(\nu)=\sum_{i=1}^M R_i^*\circ\omega_i(\nu).$$
}

\proof
Because $\Phi$ is a bijection, it is sufficient to prove that 
$$\Phi\big(\mathcal{H}(\nu)\big)=\Phi\Big(\dd\sum_{i=1}^M R_i^*\circ\omega_i(\nu)\Big)$$
for any $\nu\in cabv(T,X)$. This means to  show that, for any $f\in\mathcal{C}(T,X)$:
\begin{equation}
\tag{5.2}\label{eq:5.2}
\Phi\big(\mathcal{H}(\nu)\big)(f)=\Phi\Big(\sum_{i=1}^n R_i^*\circ\omega_i(\nu)\Big)(f).
\end{equation}

So, let us take  $f\in\mathcal{C}(T,X)$. We have successively,  using (\ref{eq:4.4}),(\ref{eq:4.5}) and the  preceding facts:
\begin{align*}
\Phi\big(\mathcal{H}(\nu)\big)(f)&=(\Phi\circ\mathcal{H})(\nu)(f)=(H_C'\circ\Phi)(\nu)(f)=\\
&=H_C'\big(\Phi(\nu)\Big)(f)=\big(\Phi(\nu)\circ H_C)(f)=\Phi(\nu)\big(H_C(f)\big)=\\
&=\int H_C(f)d\nu=\int\sum_{i=1}^M(R_i\circ f\circ\omega_i)d\nu=\int\sum_{i=1}^M(R_i\circ f)\circ\omega_id\nu=\\
&=\int\sum_{i=1}^M(R_i\circ f)d\big(\omega_i(\nu)\big)=\int\sum_{i=1}^M fd\big(R_i^*\circ\omega_i(\nu)\big)=\\
&=\int fd\Big(\sum_{i=1}^M R_i^*\circ\omega_i(\nu)\Big)=\Phi\Big(\sum_{i=1}^M R_i^*\circ\omega_i(\nu)\Big)(f),
\end{align*}

which is (\ref{eq:5.2}). \eop

C2. \underline{The countable case}

In the general framework schema, we take $\Theta=\n^*,\Sigma=\mathcal{P}(\Theta)$ and $W=card$. Consequently, we accept that
$$\sum_{i=1}^\infty\|R_i\|_o(1+r_i)<\infty.$$

\textbf{Theorem 5.2}
\textit{In the context from above, one has:
\begin{enumerate}[1.,noitemsep, nosep]
\item For any $f\in\mathcal{C}(T,X)$:
$$H_C(f)=\sum_{i=1}^\infty R_i\circ f\circ\omega_i$$
(absolute convergence in $\mathcal{C}(T,X)$).
\item For any $\nu\in cabv(T,X):$
$$\mathcal{H}(\nu)=\sum_{i=1}^\infty R_i^*\circ\omega_i(\nu)$$
(convergence in $cabv(T,X)$ with the usual variational norm).
\end{enumerate}
}

\mm
\proof
1. Due  to the inequality $\|R_i\circ f\circ\omega_i\|_\infty\leq\|R_i\|_o\cdot\|f\|_\infty$ it follows that the series $\dd\sum_{i=1}^\infty R_i\circ f\circ\omega_i$ converges absolutely in $\mathcal{C}(T,X)$.

For any $t\in T$, one has $\|(R_i\circ f\circ\omega_i)(t)\|\leq\|R_i\circ f\circ\omega_i\|_\infty$ for any $i$, hence the series $\dd\sum_{i=1}^\infty\|(R_i\circ f\circ\omega_i)(t)\|$ converges.

This shows that the Bochner integral giving $H_C(f)(t)$ is exactly
$$H_C(f)(t)=\int(R_\theta\circ f\circ\omega_\theta)(t)dW(\theta)=\sum_{i=1}^\infty(R_i\circ f\circ\omega_i)(t)$$

We can compute the sum of the absolutely convergent series (hence uniformly convergent series) $\dd\sum_{i=1}^\infty R_i\circ f\circ\omega_i$ in $\mathcal{C}(T,X)$ and now we see that
$$H_C(f)=\sum_{i=1}^\infty R_i\circ f\circ\omega_i.$$

2. The series $\dd\sum_{i=1}^\infty R_i^*\circ\omega_i(\nu)$ converges absolutely in $cabv(T,X)$, for any $\nu\in cabv(T,X)$. Indeed, for any $i$, one has
$$\|R_i^*\circ\omega_i(\nu)\|\leq\|R_i^*\|_o\cdot\|\omega_i(\nu)\|=\|R_i\|_o\cdot\|\omega_i(\nu)\|$$
(this inequality is valid, computing the respective sums on each partition of $T$). Consequently,
\begin{align*}
\sum_{i=1}^\infty\|R_i^*\circ\omega_i(\nu)\|&\leq\sum_{i=1}^\infty\|R_i\|_o\cdot\|\omega_i(\nu)\|\leq\\
&\leq\sum_{i=1}^\infty\|R_i\|_o\cdot\|\nu\|<\infty\quad a.s.o.
\end{align*}

Now, let us return to the very proof. Again it will be sufficient to prove that, for any $f\in\mathcal{C}(T,X)$ and any $\nu\in cabv(T,X)$, one has
\begin{equation}
\tag{5.3}\label{eq:5.3}
\Phi\big(\mathcal{H}(\nu)\big)(f)=\Phi\Big(\sum_{i=1}^\infty R_i^*\circ\omega_i(\nu)\Big)(f).
\end{equation}

Indeed, for such $f$ and $\nu$ (use Theorem 4.3, (\ref{eq:4.4}), (\ref{eq:4.5}) and previous remarks):

\begin{align*}
\Phi\big(\mathcal{H}(\nu)\big)(f)&=(\Phi\circ\mathcal{H})(\nu)(f)=(H_C'\circ\Phi)(\nu)(f)=\\
&=\big(\Phi(\nu)\circ H_C\big)(f)=\Phi(\nu)\big(H_C(f)\big)=\Phi(\nu)\Big(\sum_{i=1}^\infty R_i\circ f\circ\omega_i\Big)=\\
&=\sum_{i=1}^\infty\Phi(\nu)(R_i\circ f\circ\omega_i),\mbox{ using point 1.}
\end{align*}
\begin{align*}
\sum_{i=1}^\infty\Phi(\nu)(R_i\circ f\circ\omega_i)&=\sum_{i=1}^\infty\int(R_i\circ f\circ\omega_i)d\nu=\\
&=\sum_{i=1}^\infty\int(R_i\circ f)\circ\omega_i d\nu=\sum_{i=1}^\infty\int(R_i\circ f)d\big(\omega_i(\nu)\big)=\\
&=\sum_{i=1}^\infty\int f d\big(R_i^*\circ\omega_i(\nu)\big)=\int fd\Big(\sum_{i=1}^\infty R_i^*\circ\omega_i(\nu)\Big).
\end{align*}

The last equality is valid because the  series
$\dd\sum_{i=1}^\infty R_i^*\circ\omega_i(\nu)$ converges (absolutely)  in $cabv(T,X)$
and the sesquilinear uniform integral is a continuous sesquilinear map.

The last value, i.e. $\dd\int fd\Big(\sum_{i=1}^\infty R_i^*\circ\omega_i(\nu)\Big)$ is  exactly $\dd\Phi\Big(\sum_{i=1}^\infty R_i^*\circ\omega_i(\nu)\Big)(f)$ and this proves (\ref{eq:5.3}). \eop

\mm
\textbf{Remark} 
\textit{Theorems 5.1 and 5.2 show that, in the discrete case, the points of view of the present paper and of \cite{CIMN3} are dual. Namely, in the present paper, starting with the operator $H_C$, acting on continuous functions, we obtain the operator $\mathcal{H}$, acting on measures. The expression of $\mathcal{H}$ is 
$$\mathcal{H}(\nu)=\sum_i R_i^*\circ\omega_i(\nu)$$
(the sum being finite or countable).
}

\textit{Dually, in \cite{CIMN3} one starts with the operator $\overline{\mathcal{H}}$(in \cite{CIMN3} 
$\overline{\mathcal{H}}$ is denoted with $H$), acting on measures, via
$$\overline{\mathcal{H}}(\nu)=\sum_i R_i\circ\omega_i(\nu).$$
}

This operator leads naturally to the operator (acting on continuous functions) defined via the correspondence $\dd f\mapsto g=\dd\sum_i R_i^*\circ  f\circ\omega_i$. Writing $g$ in the form $g=\dd\int(R_i^*\circ f\circ\omega_i)d card$, the operator acting via $f\mapsto g$ can be viewed as the dual of $H_C$.

In the next paragraph, we shall work in the discrete case with the operator $\overline{\mathcal{H}}$ instead of $\mathcal{H}$, in the spirit of \cite{CIMN3}. This is equivalent to consider the operators $R_i^*$ instead of $R_i$. This change will not affect the exemplifications in the next paragraph, because the conditions which must be fulfilled (see the forthcoming inequalities (\ref{eq:6.1'}), (\ref{eq:6.2'}), (\ref{eq:6.3'})) are the same, due to the fact that $\|R_i\|_o=\|R_i^*\|_o$. 

\section{Invariant (fractal) measures}

Considering again the general schema, we shall construct (using the operator $\mathcal{H}$) new operators on spaces of measures and we shall look for fixed points of these new operators.

We shall call these fixed points \underline{invariant (fractal) measures}. The attribute ''invariant'' is clear. The supplementary attribute ''fractal'' will be justified further (see Example 6.2 and the Remarks following it).

An informal preliminary argument leads to the idea that, sometimes, the search of fixed points has an algebraic aspect. In this respect, one can see relation (\ref{eq:5.1}), where the fixed point equations
$$\mathcal{H}(\nu)=\nu\Rightarrow\int fd\mathcal{H}(\nu)=\int fd\nu$$
(the value $\nu(T)$ is decisive), which must be valid for any $f\in\mathcal{C}(T,X)$, lead to linear systems in the discrete finite case $\Theta=\{1,2,\dots,M\},\mbox{ for finite }T$.
(various situations can appear: no fixed points, one fixed point, many fixed points).

We shall use in the sequel the contraction principle to prove the existence and uniqueness of fixed points (invariant fractal measures).

In order to have contractions (to apply the above mentioned principle), we shall consider that one of the following conditions is fulfilled:
\begin{align}
&\int\|R_\theta\|_odW(\theta)<1 \tag{6.1}\label{eq:6.1}\\
&\int\|R_\theta\|_o(1+r_\theta)dW(\theta)<1\tag{6.2}\label{eq:6.2}\\
&\int\|R_\theta\|_o\cdot r_\theta dW(\theta)<1.\tag{6.3}\label{eq:6.3}
\end{align}

In the discrete case, these conditions become:
\begin{align}
&\sum_i\|R_i\|_o<1\tag{6.1\ghil}\label{eq:6.1'}\\
&\sum_i\|R_i\|_o(1+r_i)<1\tag{6.2\ghil}\label{eq:6.2'}\\
&\sum_i\|R_i\|_o\cdot r_i<1.\tag{6.3\ghil}\label{eq:6.3'}
\end{align}
where the sum is either $\dd\sum_{i=1}^M$ (in the finite case) or $\dd\sum_{i=1}^\infty$ (in the infinite case).

We introduce the two basic schemas used in the sequel.

\underline{First schema}

One considers a non empty set $A\subset cabv(T,X)$ such that $\mathcal{H}(A)\subset A$. One can define $\mathcal{H}_1:A\rightarrow A,\mbox{ via}$
\begin{align*}
&\mathcal{H}_1(\nu)\overset{def}{=}\mathcal{H}(\nu).
\end{align*}

Then the corresponding operator norm, denoted by $\|\mathcal{H}\|_o$, has the pro\-per\-ty $\|\mathcal{H}\|_o<1$, hence $\mathcal{H}_1$ is a contraction, because, for $\mu,\nu$ in $A$, one has
$$\|\mathcal{H}_1(\mu)-\mathcal{H}_1(\nu)\|_A\leq\|\mathcal{H}\|_o\cdot\|\mu-\nu\|_A,$$
with the corresponding norm $\|\cdot\|_A$ in $A$.

\underline{Second schema}

One considers a non empty set $A\subset cabv(T,X)$ and a measure $\mu^0\in cabv(T,X)$ having the property that
$$\mathcal{H}(A)+\mu^0\overset{def}{=}\{\mathcal{H}(\mu)+\mu^0\mid\mu\in A\}\subset A.$$
One can define $\mathcal{H}_2:A\rightarrow A$, via
$$\mathcal{H}_2(\mu)\overset{def}{=}\mathcal{H}(\mu)+\mu^0.$$

Then the corresponding operator norm, denoted by $\|\mathcal{H}\|_o$, has the pro\-per\-ty $\|\mathcal{H}\|_o<1$, hence $\mathcal{H}_2$ is a contraction, because, for $\mu,\nu$ in $A$, one has
$$\|\mathcal{H}_2(\mu)-\mathcal{H}_2(\nu)\|_A=\|\mathcal{H}\|_o\cdot\|\mu-\nu\|_A,$$
with the corresponding norm $\|\cdot\|_A$ in $A$.

For both schemas, it will be necessary to check the completeness of $A$ equipped with the metric generated by the corresponding $\|\cdot\|_A$.

In the sequel, we shall introduce some theoretical and practical exemplifications of the previous schemas.

We begin with a theoretical exemplification.

\underline{Example 6.1} (according to the second schema)

We work in the context of operator semigroups.

Let $1<N<\infty$ and consider the particular case of uniformly continuous operator semigroups on a Hilbert space $X$, given as follows:
$$R_\theta=e^{-N\theta} I,\mbox{for any }\theta\in[0,\infty).$$
Because
$$\int_0^\infty\|R_\theta\|_od\theta=\int_0^\infty e^{-N\theta}d\theta=\frac1N<1,$$
condition (\ref{eq:6.1}) is fulfilled.

In order to apply the second schema, we consider a strictly positive number $a$, hence $A=B_a(X)$ is a complete metric space for the metric given by the variational norm. Take $\mu^0\in cabv(T,X)$ such that
$$\frac aN+\|\mu^0\|\leq a\Leftrightarrow\|\mu^0\|\leq a\Big(1-\frac 1N\Big).$$

Then, for any $\mu\in B_a(X)$, one has (see Theorem 4.4): $\mathcal{H}\big(B_a(X)\big)+\mu^0\subset B_a(X)$, because
\begin{align*}
\|\mathcal{H}(\mu)+\mu^0\|&\leq\|\mathcal{H}(\mu)\|+\|\mu^0\|\leq\|\mathcal{H}\|_{o,var}\cdot\|\mu\|+\|\mu^0\|\leq\\
&\leq\|\mu\|\cdot\int_0^\infty\|R_\theta\|_od\theta+\|\mu^0\|\leq\frac aN+\|\mu^0\|\leq a.
\end{align*}

We are in position to define the contraction
$\mathcal{H}_2:B_a(X)\rightarrow B_a(X)$, given via $\mathcal{H}_2(\mu)=\mathcal{H}(\mu)+\mu^0$, for any $\mu\in B_a(X)$. Namely, for any $\mu,\nu$ in $B_a(X)$, one has
$$\|\mathcal{H}_2(\mu)-\mathcal{H}_2(\nu)\|\leq\frac 1N\|\mu-\nu\|.$$

The contraction principle says that there exists a unique fixed point $\mu^*\in B_a(X)$ of $\mathcal{H}_2$:
$$\mu^*=\mathcal{H}_2(\mu^*)=\mathcal{H}(\mu^*)+\mu^0.$$
In case $\mu^0=0$, one has $\mu^*=0$. \eop

\mm
The following three examples will refer to the discrete finite case. Namely, we shall take for the concrete illustration:
$T=[0,1],M=2$ (i.e. $\Theta=\{1,2\}$) and $\omega_1,\omega_2:[0,1]\rightarrow [0,1]$ the Cantor contractions
\begin{align*}
\omega_1(t)&=\frac t3\; \Big(r_1=\frac 13\Big)\\
\omega_2(t)&=\frac 23+\frac t3\; \Big(r_2=\frac 13\Big).
\end{align*}

It is seen that, for any $\emptyset\neq B\in\mathcal{B}$, one has
\begin{align*}
\omega_1^{-1}(B)&=(3B)\cap[0,1]\overset{def}{=}\{3t\mid t\in B\}\cap [0,1]\\
\omega_2^{-1}(B)&=(3B-2)\cap[0,1]\overset{def}{=}\{3t-2\mid t\in B\}\cap [0,1]
\end{align*}

As we said, we shall work with $\overline{\mathcal{H}}$ instead of $\mathcal{H}$ (also in the schemas' constructions) and with $R_i^*$ instead of $R_i,i=1,2$.

Each of the following three examples will be introduced theoretically, in the spirit of the aforementioned schemas and will be illustrated concretely. Proofs and computations will be merely sketched, the details being contained in \cite{CIMN3}.

\underline{Example 6.2} (according to the first schema)

Consider $X=K^n,n\in\n^*$. The hypotheses are:
\begin{enumerate}[a),noitemsep,nosep]
\item $\dd\sum_{i=1}^n R_i=I$
\item $c\overset{def}{=}\dd\sum_{i=1}^M\|R_i\|_o\cdot r_i<1$ (see (\ref{eq:6.3'})).
This is true if all $\omega_i$ are contractions and $\dd\sum_{i=1}^M\|R_i\|_o=1$.
\item $0<a<\infty$ and $v\in K^n$ is such that $\|v\|\leq a$. (hence $B_a(K^n,v)\neq\emptyset$).
\item $\emptyset\neq A\subset B_a(K^n,v)$ is such that $\overline{\mathcal{H}}(A)\subset A$ and $A$ is weak* closed. In the particular case when $\|\overline{\mathcal{H}}(\mu)\|\leq\|\mu\|$ for any $\mu\in cabv(T,K^n)$, one can take $A=B_a(K^n,v)$ (more particular, if $\dd\sum_{i=1}^M\|R_i\|_o=1$, it follows that $\|\overline{\mathcal{H}}(\mu)\|\leq\|\mu\|$ for any $\mu\in cabv(T,K^n$)). 
\end{enumerate}

Under these hypotheses, we can define $\mathcal{H}_1: A\rightarrow A$ via $\mathcal{H}_1(\mu)=\overline{\mathcal{H}}(\mu)$, for any $\mu\in A$. It follows that $\mathcal{H}_1$ is a contraction with contraction factor $\leq c$, if $A$ is equipped with the metric $d_{MK}^*$ given via $d_{MK}^*(\mu,\nu)=\|\mu-\nu\|_{MK}^*$. (according to Theorem 4.7)

Consequently, there exists a unique invariant (fractal) measure $\mu^*\in A$ of $\mathcal{H}_1$, i.e. $\mathcal{H}_1(\mu^*)=\mu^*$.

\underline{Sketch of proof}
The first basic idea is that $B_a(K^n,v)$ is a non empty compact space for the metric $d_{MK}^*$, hence $A$ is also compact for this metric, being weak* closed. The second basic idea is that condition a) guarantees the fact that $\overline{\mathcal{H}}(cabv(K^n,v))\subset cabv(K^n,v)$. Computing details and the contraction principle complete the proof.
\eop

\mm
\underline{Remarks}
\begin{enumerate}[1.,noitemsep,nosep]
\item Condition a) implies that $1=\|I\|_o\leq\dd\sum_{i=1}^M\|R_i\|_o$, hence condition \newline $\dd\sum_{i=1}^M\|R_i\|_o=1$ is extremal.

There exist situations when all the particular conditions are fulfilled (see the following Remark).
\item The classical model, producing the invariant (fractal) probability, is a particular case of Theorem 6.2 where all the particular conditions are fulfilled.

Namely, in the classical model, one has $n=1$ (hence $X=K$), $R_i\in\mathcal{L}(K)$ are given via $R_i(t)=p_it$, where all $p_i>0$ and $\dd\sum_{i=1}^M p_i=1$, hence $\dd\overline{\mathcal{H}}(\mu)=\sum_{i=1}^Mp_i\omega_i(\mu)$ for any $\mu\in cabv(T,K)$. Also, one takes $a=1,v=1$ and $A=\{\mu\in B_1(K,1)\mid\mu\geq 0\}=$ the set of all probabilities ${\mu:\mathcal{B}\rightarrow [0,1]}$. Then $A$ is weak* closed and $\dd\sum_{i=1}^M R_i=I,\sum_{i=1}^M\|R_i\|_o=\sum_{i=1}^Mp_i=1$ (hence $c<1$). For any contractions $\omega_i:T\rightarrow T,i=1,2,\dots,M$, we find a unique probability $\mu^*:\mathcal{B}\rightarrow[0,1]$ (the invariant fractal measure) such that $\dd\mu^*=\sum_{i=1}^M p_i\omega_i(\mu)$.
\end{enumerate}

\underline{Concrete illustration}
Take $n=2$ (hence $X=K^2$) and $R_1,R_2$ in $\mathcal{L}(K^2)$ such that
\begin{align*}
R_1\equiv\left(
\begin{aligned}
\alpha &\;& 0\\
0 &\;& \alpha
\end{aligned}
\right),\;
R_2\equiv\left(
\begin{aligned}
1&-\alpha&\;&0\\
&0\;&1&-\alpha
\end{aligned}
\right),
\end{align*}
where $0<\alpha<1$. Then $R_1+R_2=I,\|R_1\|_o=\alpha,\|R_2\|_o=1-\alpha$, hence $\|R_1\|_o+\|R_2\|_o=1$.

Also take $a=\sqrt{2}$ and $v=(1,1)$, hence $\|v\|=a$.

We get the invariant (fractal) measure $\mu^*=(\mu_1^*,\mu_2^*)$. Namely, the invariance equation $\mathcal{H}_1(\mu^*)=\mu^*$, i.e.
\begin{align*}
R_1\circ\omega_1(\mu^*)+R_2\circ\omega_2(\mu^*)=\mu^*\mbox{ is (for any }B\in\mathcal{B}):\\
R_1\Big(\mu^*\big((3B)\cap[0,1]\big)\Big)+R_2\Big(\mu^*\big((3B-2)\cap[0,1]\big)\Big)=\mu^*(B)
\end{align*}
In matricial form
\begin{align*}
\left(
\begin{aligned}
\alpha &\;& 0\\
0 &\;& \alpha
\end{aligned}
\right)
\left(
\begin{aligned}
\mu_1^*\big((3B)\cap[0,1]\big)\\
\mu_2^*\big((3B)\cap[0,1]\big)
\end{aligned}
\right)&+
\left(
\begin{aligned}
1&-\alpha&\;&0\\
&0\;&1&-\alpha
\end{aligned}
\right)
\left(
\begin{aligned}
\mu_1^*\big((3B-2)\cap[0,1]\big)\\
\mu_2^*\big((3B-2)\cap[0,1]\big)
\end{aligned}
\right)
=\\
&=
\left(\begin{aligned}
\mu_1^*(B)\\
\mu_2^*(B)
\end{aligned}
\right)
\end{align*}
giving, for $i=1,2$
$$\alpha\mu_i^*\big((3B)\cap[0,1]\big)+(1-\alpha)\mu_i^*\big((3B-2)\cap[0,1]\big)=\mu_i^*(B).$$

Hence $\mu_1^*=\mu_2^*=\mu$, where $\mu:\mathcal{B}\rightarrow [0,1]$ is the unique invariant (fractal) probability obtained in the classical model for $p_1=\alpha$ and $p_2=1-\alpha$. \eop

\mm
\underline{Example 6.3} (according to the second schema)

Consider $X=K^n,n\in\n^*$. The hypotheses are:
\begin{enumerate}[a),noitemsep,nosep]
\item $\dd d\overset{def}{=}\sum_{i=1}^M\|R_i\|_o(1+r_i)<1$ (see (\ref{eq:6.2'})).

This is true if all $\omega_i$ are contractions and $\dd\sum_{i=1}^M\|R_i\|_o\leq\frac 12$.
\item $0<a<\infty,\mu^0\in cabv(K^n),\emptyset\neq A\subset B_a(K^n)$ is weak* closed and one has $\overline{\mathcal{H}}(\mu)+\mu^0\in A$ for any $\mu\in A$.
In particular, if $\dd\|\mu^0\|+a\Big(\sum_{i=1}^M\|R_i\|_o\Big)\leq a$, then one can take $A=B_a(K^n)$.
\end{enumerate}

Under these hypotheses, we define $\mathcal{H}_2:A\rightarrow A$ via $\mathcal{H}_2(\mu)=\overline{\mathcal{H}}(\mu)+\mu^0$ for any $\mu\in A$. It follows that $\mathcal{H}_2$ is a contraction with contraction factor $\leq d$, if $A$ is equipped with the metric $d_{MK}$, given via $d_{MK}(\mu,\nu)=\|\mu-\nu\|_{MK}$ (according to Theorem 4.5)
Then:
\begin{enumerate}[i),noitemsep,nosep]
\item If $\mu^0=0$, it follows that $0\in A$.
\item There exists a unique invariant (fractal) measure $\mu^*\in A$ of $\mathcal{H}_2$, i.e. $\mathcal{H}_2(\mu^*)=\mu^*$. In case $\mu^0=0$, we have $\mu^*=0$.
\end{enumerate}

\underline{Sketch of proof}
Again $B_a(K^n)$ and $A$ are compact metric spaces for the metric $d_{MK}$. Computing details and the contraction principle complete the general proof. As for the particular case $\mu^0=0$, one sees that repeated application of $\mathcal{H}_2=\overline{\mathcal{H}}$ gives $\underset{n}{lim}\;\overline{\mathcal{H}}^n(\mu)=0$ for $\mu\in A$, hence $0\in A$ which is closed for $d_{MK}$. We consider also the uniqueness of $\mu^*$. \eop

\mm
\underline{Concrete illustration}

Consider $X=K^2$ and let $\mu^0:\mathcal{B}\rightarrow K^2$ act via
$$\mu^0(B)=\Big(\frac 14\lambda(B),\frac 14\delta_0(B)\Big)\mbox{ for any }B\in\mathcal{B}$$
($\lambda$ is the Lebesgue measure on $\mathcal{B}$ and $\delta_0$ is the Dirac measure concentrated at 0). Take $R_1,R_2$ in $\mathcal{L}(K^2)$ as follows:
$\dd R_i=\frac 1{10}P_i,i=1,2$, where
\begin{align*}
P_1\equiv\left(
\begin{aligned}
1 &\;& 0\\
2 &\;& 1
\end{aligned}
\right)
\mbox{ and }
P_2\equiv\left(
\begin{aligned}
1 &\;& 0\\
2 &\;& -1
\end{aligned}
\right).
\end{align*}

Consequently $\|P_1\|_o=\|P_2\|_o=1+\sqrt{2}$, giving $\dd\|R_1\|_o+\|R_2\|_o=\frac{1+\sqrt{2}}{5}<\frac 12$, so $d<1$.

Take $a=1$, hence $\dd\|\mu^0\|+a\Big(\sum_{i=1}^2\|R_i\|_o\Big)=\frac 12+\frac{1+\sqrt{2}}{5}<1=a$, because $\|\mu^0\|=\frac 12$.

The preceding theory proves the existence and uniqueness of the invariant (fractal) measure $\mu^*=(\mu_1^*,\mu_2^*)\in cabv(K^2)$.

The invariance equation is, for any $B\in\mathcal{B}$:
$$R_1\Big(\mu^*\big((3B)\cap[0,1]\big)\Big)+R_2\Big(\mu^*\big((3B-2)\cap[0,1]\big)\Big)+\mu^0(B)=\mu^*(B).$$
In matricial form
\begin{align*}
\left(
\begin{aligned}
\frac{1}{10} &\;& 0\\
\frac{2}{10} &\;& \frac{1}{10}
\end{aligned}
\right)
\left(
\begin{aligned}
\mu_1^*\big((3B)\cap[0,1]\big)\\
\mu_2^*\big((3B)\cap[0,1]\big)
\end{aligned}
\right)&+
\left(
\begin{aligned}
\frac{1}{10}&\;&0\;\\
\frac{2}{10}&\;&-\frac{1}{10}
\end{aligned}
\right)
\left(
\begin{aligned}
\mu_1^*\big((3B-2)\cap[0,1]\big)\\
\mu_2^*\big((3B-2)\cap[0,1]\big)
\end{aligned}
\right)
+\\
&+
\left(\begin{aligned}
\frac14\lambda(B)\\
\frac14\delta_0(B)
\end{aligned}
\right)=
\left(\begin{aligned}
\mu_1^*(B)\\
\mu_2^*(B)
\end{aligned}
\right)
\end{align*}
giving,
\begin{align*}
&\frac{1}{10}\mu_1^*\big((3B)\cap[0,1]\big)+\frac{1}{10}\mu_1^*\big((3B-2)\cap[0,1]\big)+\frac{1}{4}\lambda(B)=\mu_1^*(B)\\
&\frac{2}{10}\mu_1^*\big((3B)\cap[0,1]\big)+\frac{1}{10}\mu_2^*\big((3B)\cap[0,1]\big)+\frac{2}{10}\mu_1^*\big((3B-2)\cap[0,1]\big)-\\
&-\frac{1}{10}\mu_2^*\big((3B-2)\cap[0,1]\big)+\frac{1}{4}\delta_0(B)=\mu_2^*(B).
\end{align*}

Examples of computation:

$\dd\mu^*\big([0,1]\big)=\Big(\frac 5{16},\frac 38\Big),\mu^*\big(\{0\}\big)=\Big(0,\frac 5{18}\Big),\mu^*\big(\{1\}\big)=(0,0),\mu^*\Big(\Big\{\frac 23\Big\}\Big)=\Big(0,-\frac 1{36}\Big)$. \eop

\mm
\underline{Example 6.4} (according to the second schema)

We work in an arbitrary Hilbert space $X$ and consider $cabv(T,X)$ with the variational norm. Take $\mu^0\in cabv(T,X)$. The hypotheses are:
\begin{enumerate}[a),noitemsep,nosep]
\item $\dd e\overset{def}{=}\sum_{i=1}^M\|R_i\|_o<1$ (see (\ref{eq:6.1'}))
\item $\emptyset\neq A\subset cabv(T,X)$ is a closed set such that $\overline{\mathcal{H}}(\mu)+\mu^0\in A$ for any $\mu\in A$.

This is true if: either $A=cabv(T,X)$, or $A=B_a(X)$, where $0<a<\infty$ is such that $\dd\|\mu^0\|+a\Big(\sum_{i=1}^M\|R_i\|_o\Big)\leq a$.
\end{enumerate}

Under these hypotheses, we define $\mathcal{H}_2:A\rightarrow A$ via $\mathcal{H}_2(\mu)=\overline{\mathcal{H}}(\mu)+\mu^0$ for any $\mu\in A$. It follows that $\mathcal{H}_2$ is a contraction with contraction factor $\leq e$. Then:
\begin{enumerate}[i),noitemsep,nosep]
\item If $\mu^0=0$, then $0\in A$.
\item There exists a unique invariant (fractal)
measure $\mu^*\in A$ of $\mathcal{H}_2$, i.e. $\mathcal{H}_2(\mu^*)=\mu^*$. In case $\mu^0=0$, we have $\mu^*=0$.
\end{enumerate}

Similar (simpler) sketch of proof as that one given for Example 6.3.

\underline{Concrete illustration}

We begin with initial facts (here $\lambda$ is the Lebesgue measure on $\mathcal{B}$).

Any continuous function $F:[0,1]^2\rightarrow K$
(write $Q\overset{def}{=}sup\{|F(x,y)|\mid (x,y)\in[0,1]^2\}$) generates
$R:L^2(\lambda)\rightarrow L^2(\lambda)$ given via $R(\tilde{f})=\tilde{g}$, where $g:[0,1]\rightarrow K$ is the continuous function acting as follows:
$$g(x)=\int_0^1 F(x,y)f(y)d\lambda(y)\mbox{ (we work with a representative }f\in\tilde{f}).$$

Because $\|g\|_2\leq Q\|f\|_2$, we see that $R\in\mathcal{L}(L^2(\lambda))$ and $\|R\|_o\leq Q$.

Now, we shall introduce our concrete example.
Take $0<a<\infty$, $X=L^2(\lambda)$ and $F_i:[0,1]^2\rightarrow K$ continuous functions, with $Q_i\overset{def}{=}sup\{|F_i(x,y)|\mid (x,y)\in[0,1]^2\}$ and we shall assume that $\dd Q_i\leq\frac 14,i=1,2$. As previously, we shall generate, using $F_i$, the linear and continuous operators $R_i\in\mathcal{L}^2(X)$, hence $\dd\|R_i\|_o\leq Q_i\leq\frac 14,i=1,2$. Then $\dd\|R_1\|_o+\|R_2\|_o\leq\frac 12<1$.

Take also $\mu^0\in cabv(T,L^2(\lambda))$ with $\dd\|\mu^0\|\leq\frac a2$. Hence $\|\mu^0\|+a(\|R_1\|_o+\|R_2\|_o)\leq a$. 

We can apply the previous result.
The effective computation will be done for the following $F_1,F_2$ and $\mu^0$. Take $\dd F_1(x,y)=\frac{xy}{4}$ and $\dd F_2(x,y)=\frac{x^2y^2}{4}$, hence $\dd Q_1=Q_2=\frac 14$.

In order to introduce $\mu^0$, we first consider the measure $m\in cabv(T,L^2(\lambda))$ given, for any $B\in\mathcal{B}$, via $m(B)=\tilde{h}_B$ where $h_B:[0,1]\rightarrow K$ is  the continuous function acting as follows:

$$h_B(t)=\lambda(B\cap[0,t]),\mbox{ for any } t\in[0,1].$$

Then  $\dd\|m\|=\frac 23$ (see \cite{CIMN}). Finally, we take $\dd\mu^0\overset{def}{=}\frac 12m$,  hence $\dd\|\mu^0\|=\frac 13$ and $a=1$. It follows that

$$\|\mu^0\|+a(\|R_1\|_o+\|R_2\|_o)<a.$$

Hence, we obtain the unique invariant (fractal) measure $\mu^*\in cabv(T,L^2(\lambda)).$

The invariance equation is (for any $B\in\mathcal{B}$): 
$$R_1\Big(\mu^*\big((3B)\cap[0,1]\big)\Big)+R_2\Big(\mu^*\big((3B-2)\cap[0,1]\big)\Big)+\mu^0(B)=\mu^*(B).$$

In order to examine this equation, we shall consider, for any $B\in\mathcal{B}$, a representative $\tilde{f}_B\in\mu^*(B)$, thus obtaining representatives of  $R_1(\mu^*(B))$ and  $R_2(\mu^*(B))$ via (we write abusively identifying classes with representatives):
$$R_1(\mu^*(B))=\frac 14\int_0^1 xy f_B(y)d\lambda(y),\;R_2(\mu^*(B))=\frac 14\int_0^1 x^2y^2 f_B(y)d\lambda(y)$$

The invariance  equation becomes

\begin{align*}
\frac 14\int_0^1 xy f_{(3B)\cap[0,1]}(y)d\lambda(y)&+\frac 14\int_0^1 x^2y^2 f_{(3B-2)\cap[0,1]}(y)d\lambda(y)+\\
&+\frac 12\lambda(B\cap[0,x])=f_B(x),
\end{align*}
for any $B\in\mathcal{B}$ and $\lambda$-almost all $x\in[0,1]$.

In particular , taking $B=[0,1]$ and writing $f_{[0,1]}=\varphi$, we have, for almost all $x\in[0,1]$, the integral equation
$$\varphi(x)=\frac 12x+\frac 14\Big(x\int_0^1 y\varphi(y)d\lambda(y)+x^2\int_0^1 y^2\varphi(y)d\lambda(y)\Big)$$
with the solution
$$\varphi(x)=\frac{24}{3329}(76x+5x^2).$$
\eop

\mm
Similar results to those introduced in the last three examples, but in the countable discrete case, can be obtained under conditions of type 
(\ref{eq:6.1'}), (\ref{eq:6.2'}) and (\ref{eq:6.3'}) with infinite sums. (see \cite{CIMN3}).

We shall introduce a result from \cite{CIMN3} to illustrate this point of view.

For a general Hilbert space, consider an arbitrary  $P\in\mathcal{L}(X)$ and define the sequence $(R_i)_{i\geq 1}\subset\mathcal{L}(X)$ via $\dd R_i=-\frac 1{i!}P^i$ (here $P^i=P\circ P\circ\dots\circ P\; i$ times).

Let $\mu^0\in cabv(T,X)$. Assume that $(t_i)_{i\geq 1}$ is a sequence in $T$ with $t_i$ distinct such that all $\omega_i$ are constant ($\omega_i(t)=t_i$ for any $t\in T$ and any $i$).

Then one can see that (according to the second schema) the formula of $\mathcal{H}_2:cabv(T,X)\rightarrow cabv(T,X)$ is
$$\dd\mathcal{H}_2(\mu)=-\sum_{i=1}^\infty\frac 1{i!}P^i\big(\mu(T)\big)+\mu^0$$ and $\mathcal{H}_2$ possesses the fixed point $\mu^*\in cabv(T,X)$ given  via
$$\mu^*=-\sum_{i=1}^\infty\frac 1{i!}\delta_{t_i}\big(P^i\circ exp(-P)\big)\big(\mu^0(T)\big)+\mu^0.$$

\Addresses

\end{document}